\documentclass[15pt]{amsart}
\usepackage{graphicx} 
\DeclareGraphicsExtensions{.pdf,.eps,.fig}
\usepackage{url}
\usepackage{amsmath}
\usepackage{amssymb}
\usepackage{graphicx}
\usepackage{subcaption}
\usepackage{bm}
\usepackage{amsthm}
\usepackage{amscd}
\usepackage{amsfonts}
\usepackage{adjustbox}
\usepackage{fancyhdr}
\newcommand{\innerproduct}[2]{\langle #1, #2 \rangle}

\usepackage{hyperref}
\usepackage{float}
\usepackage{tikz}
\usetikzlibrary{mindmap,trees}

\DeclareGraphicsExtensions{.pdf,.eps,.fig}

\usepackage{graphicx} 
\DeclareGraphicsExtensions{.pdf,.eps,.fig,.pdf_tex}
\usepackage{color}
\usepackage{amsmath}
\usepackage{amssymb}
\usepackage{amsthm}
\usepackage{amscd}
\usepackage{amsfonts}
\usepackage{fancyhdr}
 \usepackage{hyperref}
 \usepackage{float}
 \usepackage[all]{xy}
 \usepackage{tikz-cd}
\usepackage{mathtools}
\usepackage{amsfonts}
\usepackage{enumerate}
\usepackage{enumitem}
\usepackage{comment}
\usepackage{leftindex}

\newcommand{\derfunc}[2]{ \frac{\delta #1 }{\delta #2}}
\newcommand{\partiald}[2]{\frac{\partial #1 }{\partial #2}}
\newcommand{\norm}[1]{\lVert #1  \rVert}
\theoremstyle{plain} \numberwithin{equation}{section}
\newtheorem{theorem}{Theorem}[section]

\newtheorem{proposition}[theorem]{Proposition}
\theoremstyle{definition}
\newtheorem{definition}[theorem]{Definition}

\newtheorem{remark}[theorem]{Remark}
\newtheorem{example}[theorem]{Example}
\newcommand{\blue}{\color{blue}}
\newcommand{\black}{ \color{black} }

\title{Symplectic Geometry, Poisson Geometry, and Beyond}
\author{Ivan Contreras, Diego Martinez, Nicolas Martinez, Diego Rodriguez}

\begin{document}
\maketitle
\begin{center}
\emph{``Everything is a Lagrangian submanifold". Alan Weinstein}
\end{center}

\begin{abstract}
    Symplectic and Poisson geometry emerged as a tool to understand the mathematical structure behind classical mechanics. However, due to its huge development over the past century, it has become an independent field of research in differential geometry. In this lecture notes, we will introduce the essential objects and techniques in symplectic geometry (e.g  Darboux coordinates, Lagrangian submanifolds, cotangent bundles) and Poisson geometry (e.g symplectic foliations, some examples of Poisson structures). This geometric approach will be motivated by examples from classical physics, and at the end we will explore applications of symplectic and Poisson geometry to Lie theory and other fields of mathematical physics.
\end{abstract}

\section{Introduction}
These notes are based on the four lectures of the mini-course taught by the first and third author at Universidad Nacional de Colombia (July 22-26, 2024) during the first ECOGyT (``Encuentro Colombiano de Geometr\'ia y Topolog\'ia"), and by the first author at University of Notre Dame (June 17-21, 2024), during the Undergraduate Workshop as part of the Thematic Program in Field Theory and Topology. The audience for these two events was quite diverse, ranging from undergraduate math students, and graduate students in different fields besides geometry, to physics students. For this reason, the lectures were intended to introduce the field with a small glimpse to recent results, and to motivate a deeper understanding for those interested in the field of mathematical physics. The first step of the lectures is to motivate and to define the Poisson bracket associated with the dynamics of Hamiltonian systems. Later on we will see that this bracket also encodes the geometric structure of the so called Poisson manifolds, the main object of study in Poisson geometry. After studying the main features of Poisson geometry, symplectic manifolds will appear as a natural example, since all symplectic manifolds are Poisson. As an illustrative example, the phase space of a physical system (the space with physical momentum and position as coordinates) is a symplectic manifold, and it is equipped with a natural Poisson bracket that is used to study of the equations of motion. This manuscript is not intended to cover a whole course on Poisson geometry or Hamiltonian systems, but we aim to introduce the basics of the subject to motivate the reader to delve deeper into this area and its related topics on mathematics or physics (or both).

Symplectic and Poisson geometry can be extended beyond physics and geometry. There are numerous connections with other fields of mathematics, and some of such connections are depicted in the diagram below, which is not as complete as it is desirable, just because connections with other fields and areas are continuously increasing in number and complexity to be included in one single diagram.

\begin{adjustbox}{max totalsize={.9\textwidth}{.7\textheight},center}
\begin{tikzpicture}
  \path[mindmap,concept color=black,text=white]
    node[concept] {Symplectic \& Poisson geometry}
    [clockwise from=0]
    child[concept color=green!50!black] {
      node[concept] {Physics}
      [clockwise from=90]
      child { node[concept] {Field Theories: Classical \& Topological} }
      child { node[concept] {Sigma Models} }
      child{node[concept] {Poisson-Lie T-Duality}}
    }  
    child[concept color=blue] {
      node[concept] {Analysis}
      [clockwise from=-30]
      child { node[concept] {PDEs} }
      child{ node[concept] {Symplectic integrators}}
    }
    child[concept color=red] { 
      node[concept] {Algebra} 
      [clockwise from=-90]
      child { node[concept] {Lie Algebras and bialgebras} }
      child { node[concept] {Cluster algebras} }
      child{ node[concept] {D-modules }}
    }
    child[concept color=orange] { 
      node[concept] {Combinatorics}
      [clockwise from=180]
      child { node[concept] {Delzant Polytopes} }
      child { node[concept] {Toric Manifolds} }
    }
    child[concept color=violet] {
      node[concept] {Knot Theory \& Topology}
      [clockwise from=0]
      child { node[concept] {Floer Theory} }
    };
\end{tikzpicture}
\end{adjustbox}	

\subsection{Notation} For consistency, throughout this work, we will assume that all configuration spaces are finite-dimensional manifolds, except for certain special cases that will be discussed in the first section. Additionally, all functions, vector fields, and tensors are considered smooth.

\subsection{Suggested References} For more detailed information on manifolds and smooth structures, we refer the reader to \cite{lee2012smooth, jost2017riemannian, JMMS}. For more on symplectic geometry, the references \cite{Cannas,WBQuanti,LSM,ps} are highly recommended, and for more details on Poisson geometry, see \cite{PoissCFM, dufour2005poisson, laurent-gengoux2013poisson}.   

\subsection*{Acknowledgments} We thank the organizers of ECOGyT 2024 and the CMND Thematic Program in Field Theory and Topology, for creating an excellent learning environment, and collaborative work between mathematicians and physicists. IC thanks Universidad Nacional de Colombia and University of Notre Dame, and NM thanks Amherst College for the hospitality.

\section{Hamiltonian Systems}\label{secham}
Let us start our discussion with  Newton's second law of motion of a particle with normalized mass, which describes the dynamics of a classical system:
\begin{equation}
    F = \ddot{q} \footnote{Here, we are taking $m=1$}
    \label{nsl}
\end{equation}
Here $q\in \mathbb{R}^{n}$ denotes the position of some particle, and $F$ is the total force exerted on the particle. Note that equation \ref{nsl} is a second-order ordinary differential equation, which we wish to simplify  to a system of first-order ordinary differential equations. If the force $F$ is conservative \footnote{This means that the work is independent on the path and the time, and it only depends on the initial and end points} then we can write Newton's second law in terms of a physical potential $V(q)$:
\begin{equation}
    \ddot q = - \frac{\partial V(q)}{\partial q}
\label{Newton}
\end{equation}
Introducing linear momentum $p = \dot{q}$ \footnote{In the presence of a magnetic field, momentum will be affected, and this equation won't hold.  }, we can write Newton's second law as a system of two first-order differential equations:
\begin{equation}
    \begin{split}
        p = \dot{q} \\
        - \frac{\partial V(q)}{\partial q} = \dot{p}
    \end{split}
    \label{eom}
\end{equation}

The total energy, which in this case is simply the sum of kinetic and potential energy $H: \mathbb{R}^{2n} \to \mathbb{R}$, 
\begin{equation}
    H(q,p) = \frac{1}{2}\sum_{i=1}^{n}p_{i}^{2} + V(q)
\label{H(q,p)}
\end{equation}
will be called the {\it Hamiltonian function}. Note that if we take the derivatives of $H$ with respect to both $p$ and $q$, we get the following system:
\begin{equation}
\begin{split}
    \dot{q} = \frac{\partial H}{\partial p} \\
    \dot{p} = -\frac{\partial H}{\partial q}
\end{split}
\label{hje}
\end{equation}
which is known as the \emph{Hamilton equations}. To see the usefulness of having a Hamiltonian function associated with a system of Ordinary Differential Equations, notice that if we suppose that $\phi(t)=(q(t),p(t))$ is a solution, then
\begin{equation*}
    \frac{d}{dt}H(\phi(t)) = \nabla H (\phi(t))\cdot \dot{\phi}(t) = \frac{\partial H}{\partial q}\dot{q} + \frac{\partial H}{\partial p}\dot{p} = 0.     
\end{equation*}
That is, the total energy is a conserved quantity for any isolated classical system, or in other words, the solutions of the equations of motion \ref{eom} are contained in the level sets of the Hamiltonian, that is, all the points $c\in \mathbb{R}$ such that $H(q,p)= c$. 

To illustrate this property, we consider the example of a simple pendulum.

\begin{example} (Pendulum)
    The motion of a pendulum with length $L$, can be described in terms of {\black an} angular position $\theta$, the angular velocity $\omega${\black $=\dot \theta$}, both as functions of time $t$, and the acceleration due to gravity $g$:
    \begin{equation*}
        \begin{cases}
    \ddot{\theta}(t) = - \frac{g}{L} \sin(\theta(t))   \\
\dot{\theta}(t) = \omega(t) 
\end{cases}
    \end{equation*}
    The Hamiltonian function associated to this system is:
    \begin{equation*}
            H = \frac{{\black \omega}^{2}L^{2}}{2}+gL\cos{(\theta)}
    \end{equation*}
    An analytical solution of the system comes up with elliptic functions, we can also study the level curves of the Hamiltonian $H$ to have a qualitative description of the dynamics. The equation that gives the level curves of the Hamiltonian is:
    \begin{equation*}
            c = \omega^{2}+2\frac{g}{L}\cos{(\theta)}
    \end{equation*}
    \begin{figure}
        \centering
        \includegraphics[width=1\linewidth]{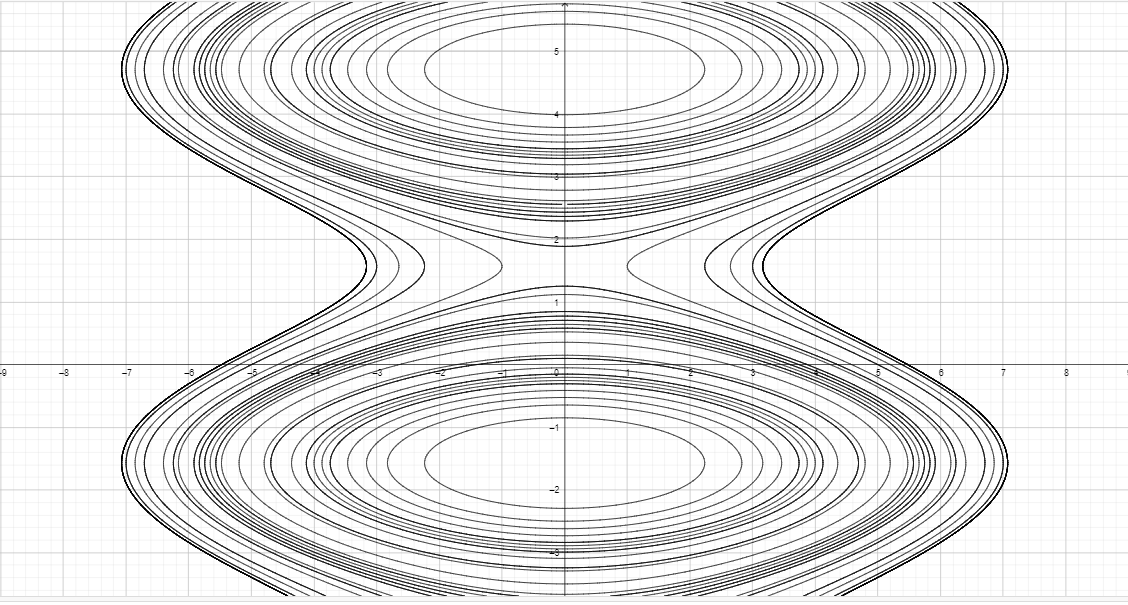}
        \caption{Phase diagram of the pendulum.}
        \label{pendpd}
    \end{figure}
    In Figure \ref{pendpd} we can see the periodic orbits that represent the low-energy oscillations, as energy increases the the orbits stop being periodic which reflects a more complex dynamics. 
$\diamond$\end{example}

Notice that we can write the system \ref{hje} using matrix notation as follows. Let $z = (q_{1}, \dots, q_{n}, p_{1}, \dots , p_{n} ) \in \mathbb{R}^{2n} $, then:
\begin{equation}
            \dot{z} = -J_{0}\nabla H(z)  = \begin{pmatrix}
    \partial H / \partial p\\
    -\partial H / \partial q
\end{pmatrix}
\label{hamjac}
        \end{equation}
 With $J_{0}$  the matrix:
            \begin{equation*}
                J_{0}=
                \begin{pmatrix}
                    0 & \mathbf{I_{n}} \\
                    -\mathbf{I_{n}} & 0
                \end{pmatrix}
            \end{equation*}
 and  $\mathbf{I_{n}}$ the $n\times n$ identity matrix. We can also write the equation \ref{hamjac} in terms of the \emph{Hamiltonian vector field}:
 \begin{equation*}
      \dot{z} = X_{H}{\black( z)}
 \end{equation*}
A geometric interpretation of $J_{0}$ can be given if we take $n=1$:
 \begin{equation*}
                J_{0}=
                \begin{pmatrix}
                    0 & 1 \\
                    -1 & 0
                \end{pmatrix}
            \end{equation*}
 Note that equation \ref{hamjac} this is a $\frac{\pi}{2}$ rotation of the gradient $\nabla H(z)$. This rotation makes the gradient point in the direction of the level curves of the Hamiltonian. 
 
 Matrix $J_{0}$ being {\black skew}-symmetric suggests the relation of classical mechanics to symplectic forms, which will be addressed in subsequent sections.

 Another important ingredient of the Hamiltonian formulation of classical mechanics is the Poisson bracket operation between functions:
 \begin{equation}
     \{f,g\} =  \sum _{i=1}^{n} \left ( \partiald{f}{q_{i}} \partiald{g}{p_{i}} - \partiald{f}  {p_{i}} \partiald{g}{q_{i}}\right)
     \label{PoissBrack}
 \end{equation}
 The evolution of any function in the system can be determined with the Hamiltonian function $H$ and the Poisson bracket $\{,\}$ with the relation:
 \begin{equation}\label{eq:dynamBRCKT}
      \dot{f} = \{f,H\}
 \end{equation}
The Poisson bracket defined above can be used to determine the dynamics of smooth functions in $\mathbb{R}^2$ (i.e $f\in C^{\infty}(\mathbb{R}\times \mathbb{R}^{2n})$), such as the Hamiltonian $H$. This is because the phase space of our initial system is $\mathbb{R}^{2n}$. If we wanted to study the equations of motion of a system on a manifold $M$ we need to understand the theory and structure behind the Poisson bracket and the Hamiltonian formulation. We will explore this approach in more detail in the next section.
Now, we will give examples of Hamiltonian systems with different bracket operations.
\begin{example}(Euler's top) \label{rigidexample}
    The motion of a top in the absence of gravity, moving around its center of mass is defined by Euler's equations:
    $$\begin{cases}
\dot{L_1}=\dfrac{I_2-I_3}{I_2I_3}L_2L_3,&\\
\dot{L_2}=\dfrac{I_3-I_1}{I_1I_3}L_1L_3,&\\
\dot{L_3}=\dfrac{I_1-I_2}{I_2I_1}L_2L_1,&\\
\end{cases}$$
Where {\black $L(t)$} $\in \mathbb{R}^{3}$ is the angular momentum and $I_{1},I_{2},I_{3}$ are the moments of inertia. This is a Hamiltonian system related to the Hamiltonian:
\begin{equation*}
    H(L)=\dfrac{1}{2}\sum_{1=1}^3\dfrac{L_i^2}{I_i}
\end{equation*}
and to the Poisson bracket operation:
\begin{equation*}
\{f,g\}(L)=\langle \nabla_L f\times \nabla_L g, L\rangle
\end{equation*}

To study rigid motion with more detail, check \cite{ABRAHAM,JMMS}.
$\diamond$\end{example}

\begin{example} (Lotka-Volterra equations)
    The Lotka-Volterra equations generalize the famous 2-species (predator-prey) dynamics to a $n$ biological species dynamics. A system of $n$ non-linear differential equations describes the model \cite{LotVol}:
    \begin{equation}
        \dot{x}_{i} = \epsilon_{i} x_{i} + \sum_{j=1}^{n}a_{ij}x_{i}x_{j}
        \label{lotka}
    \end{equation}
    Where $x_{j}$ is the population density of the species $j$, the $a_{ij}$ terms describe the interaction between species $i$ and $j$. The factor $\epsilon_{i}$ is related to the growth rate of the species $i$ due to its environment, for example, if $\epsilon_{i}>0$ then the species can survive and increase by acquiring food from the environment, but $\epsilon_{i}<0$ means that the species cannot survive and its population will decrease, clearly $\epsilon_{i}=0$ implies that the population will stay constant if $a_{ij}=0$. Considering an skew-symmetric interaction matrix $a_{ij}=-a_{ji}$ and the existence of an equilibrium $(q_{1}, \dots, q_{n})$ satisfying:
    \begin{equation*}
        \epsilon_{i} + \sum_{j=1}^{n}a_{ij}q_{j} = 0 
    \end{equation*}
    we can write the system \ref{lotka} as a Hamiltonian system, for the Hamiltonian function:
    \begin{equation*}
        h = \sum_{i}^{n}\left( x_{i} - q_{i}\log(x_{i}) \right)
    \end{equation*}
    and the bracket operation:
    \begin{equation*}
        \{f,g\} = \sum_{i<j} a_{ij}x_{i}x_{j}\left( \partiald{f}{x_{i}}\partiald{g}{x_{j}} -  \partiald{f}{x_{j}}\partiald{g}{x_{i}}  \right)
    \end{equation*}
    In this way, the Lotka-Volterra system \ref{lotka} is determined by:
    \begin{equation*}
        \dot{x}_{i} = \{x_{i},h\}
    \end{equation*}
$\diamond$\end{example}

The Hamiltonian systems are not uniquely defined for finite-dimensional configuration spaces, but we can also define them in the infinite-dimensional case. Here we present some examples but we refer to \cite{khesin_wendt_2009} for more details.

\begin{example}(Incompressible fluid) \label{example: fluid}
\label{EulerEqn}
The famous Euler equations give the physical description of a homogenous, incompressible fluid, with velocity $\textbf{v}$ and pressure $p$:
\begin{equation*}
    \partiald{\textbf{v}}{t} + (\textbf{v}\cdot \nabla)\textbf{v} = -\nabla p 
\end{equation*}
If we consider the kinetic energy of the fluid, as its Hamiltonian function:
\begin{equation*}
    H = \frac{1}{2}\int_{D} \norm{\textbf{v}}^{2} d^{3}x
\end{equation*}
With $D\subseteq \mathbb{R}^{3}$ some region in which the fluid is contained. Euler equations can be obtained as $\dot{x}=\{x,H\}$ relative to the next bracket \cite{JMMS}:
\begin{equation*}
    \{F,G\}(\textbf{v}) = \int_{D} \textbf{v} \left[\derfunc{F}{\textbf{v}} , \derfunc{G}{\textbf{v}} \right] d^{3}x
\end{equation*}
The functional derivatives can be computed as the usual Fr\'echet derivative and the bracket {\black $\left [ \cdot , \cdot \right ]$} is the Jacobi-Lie bracket for vector fields.
$\diamond$\end{example}
\begin{example}(Vl\'asov-Maxwell system)
\label{VlasovMaxwell}
The basic model to describe the physics of a plasma is the Vl\'asov-Maxwell system. A plasma is a gas heated to the point in which the thermal energy is greater than the electric bonding energy, consequently, electrons leave their atoms off. The time evolution of the system is defined by the Vl\'asov equation, in terms of a density function $f(\textbf{x},\textbf{v},t)$ and the electric and magnetic fields $\textbf{E}, \textbf{B}$:
\begin{equation*}
    \partiald{f}{t}+\textbf{v}\cdot \nabla_{x}f + \frac{q}{m}\left( \textbf{E}+ \textbf{v}\times \textbf{B}\right) \cdot \nabla_{v}f = 0
\end{equation*}
The evolution of $\textbf{E}, \textbf{B}$ is given by Maxwell equations:
\begin{align*}
    -\frac{1}{c^{2}}\partiald{\textbf{E}}{t} + \text{curl}(\textbf{B}) &= \mu_{0}\int \textbf{v}f(\textbf{x},\textbf{v},t) d\textbf{v}  \\
    \partiald{\textbf{B}}{t} +  \text{curl}(\textbf{E}) &= 0 \\
     \text{div}(\textbf{E}) &= \frac{1}{\varepsilon_{0}} \int f(\textbf{x},\textbf{v},t) d\textbf{v}\\
    \text{div}(\textbf{B}) &= 0
\end{align*}
It is possible to write Vl\'asov-Maxwell equations as a Hamiltonian system relative to a Poisson bracket \cite{VM}. Therefore, the Vl\'asov-Maxwell system is a Hamiltonian system relative to the next bracket:
\begin{multline*}
     \{F,G\}_{MV}(f,\textbf{E},\textbf{B}) =\int f \left \{\derfunc{F}{f}, \derfunc{G}{f} \right \}_{xv}  d\textbf{x}d\textbf{v}+ \int \left(\derfunc{F}{E}\cdot \text{curl}\derfunc{G}{B} - \derfunc{G}{E}\cdot \text{curl}\derfunc{F}{B} \right) d\textbf{x} \\ 
         + \int \left(\derfunc{F}{E}\cdot \partiald{f}{v}\derfunc{G}{f} -\derfunc{G}{E}\cdot \partiald{f}{v} \derfunc{F}{f}\right)d\textbf{x}d\textbf{v} + \int fB \cdot \left(\partiald{}{v}\derfunc{F}{f}\times\partiald{}{v}\derfunc{G}{f} \right)d\textbf{x}d\textbf{v}
\end{multline*}
and the Hamiltonian:
\begin{equation*}
     H(f,E,B)= \frac{1}{2}\int \norm{\textbf{v}}^{2}f(\textbf{x},\textbf{v}) d\textbf{x} d\textbf{v}+\frac{1}{2} \left( \int 	\lVert \textbf{E} \rVert^{2}  + 	\lVert \textbf{B} \rVert^{2} d\textbf{x} \right)
\end{equation*}
$\diamond$\end{example}
Therefore, Vl\'asov-Maxwell equations are $\dot{F} = \{F,H\}$.

\subsection{The Poisson bracket} \label{subsection: poisson bracket}
Along the previous examples, we identify a bracket operation defining the dynamic of the particular physical systems. However, the bracket operation also has other properties that maybe are not easily verify at a first sight, but they define the main properties of the physical system and also will allow to develop an underlying geometry. We now present one of the main definitions of these notes.
\begin{definition}
    A \textbf{Poisson manifold} is a manifold $M$ equipped with a Poisson bracket operation $\{,\}$ defined on the space of smooth functions $C^{\infty}(M)$.
     \begin{equation*}
        \{  ,\} \colon C^{\infty}(M) \times C^{\infty}(M) \mapsto C^{\infty}(M)
    \end{equation*}
      that gives $C^{\infty}(M)$ a Lie algebra structure, that satisfies Leibniz identity:
    \begin{equation*}
        \{ f , g\cdot h\} = g\cdot\{f, h\} +\{f,g\}\cdot h
    \end{equation*}
     for all $f,g,h \in C^{\infty}(M)$.
\end{definition}
The fact that $\{,\}$ gives $C^{\infty}(M)$ a Lie algebra structure means that the bracket also satisfies the next properties:
\begin{itemize}
    \item \textbf{$\mathbb{R}$ Bilinearity:} $\{\lambda f + \mu g , h\} = \lambda \{f , h \}+\mu \{g,h\}$ para $\lambda , \mu \in \mathbb{R}$
    \item \textbf{Skew-Symmetry:} $\{f,g\} = - \{g,f\}$
    \item \textbf{Jacobi's identity:} $ \{f,\{g,h\}\}+\{h,\{f,g\}\}+\{g,\{h,f\}\} = 0 $
\end{itemize}
We can have a local structure in coordinates for the bracket for any Poisson manifold $(M,\{ ,\})$. If $U$ is an open subset of $M$, and $(U, x_{1}, \dots , x_{n})$ is a local chart, then for some functions $\pi_{ij} \in C^{\infty}(M)$:
\begin{equation*}
    \{f,g\}\big|_{U} = \sum_{i,j=1}^{n}\pi_{ij} \partiald{f}{x_{i}}\partiald{g}{x_{j}}
\end{equation*}

The Leibniz identity suggests that the Poisson bracket is a derivation, a classical result from differential geometry is that there is a 1-1 correspondence between derivations and vector fields, that motivates the next definition.
\begin{definition}\label{def: ham poiss}
    Let $f \in C^{\infty}(M)$, there exists a vector field $X_{f} \in \mathfrak{X}(M)$ such that:
    \begin{equation*}
        X_{f}(g) =\{f,g\}
    \end{equation*}
    for all $g \in C^{\infty}(M)$. $X_{f}$ is called the \textbf{Hamiltonian vector field} of $f$.
\end{definition}
\begin{example}
    If we consider our Poisson manifold to be $M= \mathbb{R}^{2n}$, i.e the canonical phase space, given a function $f\in C^{\infty}(M)$, its Hamiltonian vector field $X_{f} \in \mathfrak{X}(\mathbb{R}^{2n})$ is:
    \begin{equation*}
        X_{f} = \sum _{i=1}^{n} \left ( \partiald{f}{q_{i}} \partiald{}{p_{i}} - \partiald{f}  {p_{i}} \partiald{}{q_{i}}\right)
    \end{equation*}

In particular, if we take $H \in C^{\infty}(M)$ as a Hamiltonian function, then we have that its Hamiltonian vector field will generate the time evolution of a time-dependent function $f$ in $M$, that is:
\begin{equation*}
    X_{H}(f)=  \frac{df}{dt}
\end{equation*}
$\diamond$\end{example}

Also, the hamiltonian vector field induce a Lie algebra homorphism of algebras, thanks to the Jacobi identity, between $C^{\infty}(M)$ and the set of Hamiltonian vector fields $\mathfrak{X}(M)$:
\begin{proposition}
$[X_{f},X_{g}] = X_{ \{f,g\}}$ is equivalent to a Lie algebra homomorphism between $C^{\infty}(M)$ and $\mathfrak{X}(M)$. 
\end{proposition}

\begin{definition}\label{Casimir}
    A \textbf{constant of motion}, \textbf{integral of motion} or a \textbf{first integral} is a function $f \in C^{\infty}(M)$ such that $\{H,f\}=0$, with $H$ the Hamiltonian function. That is to say, $f$ is constant along the integral curves of $X_{H}$. A \textbf{Casimir function} of the Poisson structure, is a function $C\in C^{\infty}(M)$ such that $\{C,f\}=0$, for all $f \in  C^{\infty}(M)$. Stated differently, its Hamiltonian vector field is $X_{C}=0$.

\end{definition}
Another import consequence of the Jacobi identity is he relation between two contant of motion, as the next property states:
\begin{proposition}
    If $f,g \in C^{\infty}(M)$ are constants of motion, then the function $\{f,g\}$ is also a constant of motion.
\end{proposition}

\begin{remark}
 This is one of the foundational results on classical mechanics due to Jacobi, who proved the same assertion but for the canonical bracket  in \eqref{PoissBrack}.
\end{remark}


The concept of Casimir function is related to conserved quantities in physics. To illustrate this, let's consider the following example.

\begin{example}
    In example \ref{rigidexample} we already established a Poisson structure for the Euler top. The Euler top is a rigid body, which is a body rotating through a fixed point, in such cases, the usual suspect for a conserved quantity is the angular moment, indeed $C(L)= \frac{\norm{L}^{2}}{2}$ is a Casimir function.  
$\diamond$\end{example}

\section{The background geometry for Hamiltonian systems } \label{poissonstr}

In the previous section, we have seen instances of Hamiltonian systems that are essentially different. For instance, the Euler top has non-trivial Casimir functions, while the system with potential force only has trivial (i.e constant) Casimir functions. This leads us to begin with a simplified version of the Hamiltonian system and its geometry: the case of trivial Casimir functions, that are known as symplectic Hamiltonian systems. 

In this section, we will focus on the geometric framework essential for the formal description of Hamiltonian systems. We begin by introducing the symplectic structure on a vector space, setting the foundation for its generalization to non-linear structures, namely manifolds. Following this, we explore Poisson geometry, which arises as a degenerate extension of symplectic structures. We will examine methods for constructing Poisson structures through foliations and highlight the significance of the Lie-Poisson structure in physical applications.

\subsection{Symplectic structures}
In the previous section, we expressed Hamilton's equations \ref{hje} using a skew-symmetric matrix \( J_0 \), which, as we will show, is the matrix representation of what is known as a symplectic form on a vector space. We will begin by examining the fundamental properties of symplectic structures in the context of vector spaces. Subsequently, we will introduce their non-linear counterparts on manifolds, known as symplectic manifolds, where these structures find broader and more intricate applications, such as a generalization of classical mechanics for any phase space.  

For simplicity, we refer to \cite{Cannas, mcduff_salamon, SympWeinstein} and the references therein as standard sources for this section.
\subsubsection{\textbf{Symplectic linear algebra} }\index{symplectic linear algebra}
We start by considering geometric structures on vector spaces (which may be infinite-dimensional). For the purpose of these notes, we restrict ourselves to vector spaces over $\mathbb R$. 
\begin{definition}\label{def:nondeg}\label{def: symplectic}
Let $V$ be a vector space over $\mathbb R$. A skew-symmetric form $\omega \in \Lambda^2(V^*)$ is called \textit{non-degenerate} \index{non degenerate form} if 
the induced linear map
\begin{eqnarray*}
\omega^{\#}\colon V &\to& V^*\\
\omega^{\#}(v)(w)&:=& \omega(v,w)
\end{eqnarray*}
is an isomorphism, where $V^*$ denotes the linear dual space of $V$. In this case, $\omega$ is called {\it symplectic} form and $(V,\omega)$ a symplectic vector space.
\end{definition}

A more relaxed version of non-degeneracy is required for infinite-dimensional vector spaces.
\begin{definition}\label{def: weaklynondeg}\label{weak}
A form $\omega$ is called \emph{weakly non degenerate} \index{weakly non degenerate} if $\omega^{\#}$ is an injective map. A bilinear form  $\omega$ on $V$ is called \emph{weak symplectic} \index{weak symplectic} if it is skew-symmetric and weakly non-degenerate. A vector space $V$ equipped with a weakly symplectic form $\omega$ is called a \emph{weakly symplectic vector space}\index{weakly symplectic vector space}.
\end{definition}

Observe that, in the finite-dimensional case, the facts of $\omega^{\#}$ being bijective and $\omega$ being non degenerate are equivalent. This implies for example, that the finite dimensional symplectic vector spaces are even dimensional.
\begin{example}[Standard symplectic structure]\index{Standard Symplectic Structure}\label{ex:standard_symp} Let $V=\mathbb R^{2n}$. Consider the basis 
\[\beta=\{q_1,q_2,\cdots q_n, p^1, p^2, \cdots, p^n\},\]
and the bilinear form $\omega$ defined by
\begin{eqnarray}
\omega(q_i, q_j)&=&\omega(p^i,p^j)=0, 1\leq i,j \leq n,\\
\omega(q_i, p^j)&=& \delta_{ij}, 1\leq i\leq j \leq n\\
\omega(p^j, q^i)&=& -\delta_{ij}, 1\leq j\leq i \leq n.
\end{eqnarray}
Then, $(V,\omega)$ is a symplectic vector space. On this basis, the symplectic form is represented by the matrix
\begin{equation}
    \omega \mid_{\beta}=\begin{bmatrix} 0&\mathbf{I_{n}}\\-\mathbf{I_{n}}&0
    \end{bmatrix}
\end{equation}
$\diamond$\end{example}
Notice that this matrix is the same that we use to write Hamilton equation \ref{hamjac} for canonical coordinates in $\mathbb{R}^{2n}$, for a general manifold $M$ we need to extend the idea of for non-linear structures. 

A coordinate-free way to describe the previous construction is using dual spaces:
\begin{example}[Linear canonical space]\index{Cotangent Space}\label{ex:cotangent_vect} 
Let $(W, \langle \cdot, \cdot \rangle)$ be a finite-dimensional vector space with inner product $\langle \cdot, \cdot \rangle$. Let $V=W \oplus W^*$. Given $x=(w_1, \alpha_1), y=(w_2, \alpha_2) \in V$, we define the following bilinear form:
\begin{equation}
    \omega(x,y)=\alpha_1(w_2)-\alpha_2(w_1).
\end{equation}
Then, $(V,\omega)$ is a symplectic vector space.

$\diamond$\end{example}

\begin{definition}\label{symplectomorphism} A bijective linear map $f: (V,\omega_V)\to (W, \omega_W)$ is called a \emph{symplectomorphism} \index{symplectomorphism} if $f^*\omega_W=\omega_V$\footnote{This is an extension of the transpose map, i.e $f^*\omega_W(x,y)=\omega_V(f(x),f(y))$}. We denote the set of linear symplectomorphisms between symplectic vector spaces of dimension $d$ and over $\mathbb{R}$ by Sp$(d,\mathbb{R})$.
\end{definition}

\begin{example}\label{ex:symplectomorphism} The symplectic vector spaces in Examples \ref{ex:standard_symp} and \ref{ex:cotangent_vect} are symplectomorphic.
$\diamond$\end{example}
\begin{remark}
    Weyl originally introduced the term 'complex groups' to describe the group of all linear transformations that preserve the structure of a bilinear, non-degenerate, skew-symmetric form. However, to avoid confusion with complex numbers, he later replaced the word ``complex" with its Greek equivalent, ``symplectic" \cite{weyl1946classical}.
\end{remark}

The following result shows that there is a normal form for symplectic structures on vector spaces.
\begin{theorem}[Darboux basis]\label{thm:Darboux basis} \index{Darboux basis} Let $(V,\omega)$ be a symplectic vector space. Then there is a basis 

\[(q_1,q_2\cdots q_n, p^1, p^2, \cdots, p^n)\] of $V$ such that $(V,\omega)$ is symplectomorphic to the standard symplectic vector space in Example \ref{ex:standard_symp}.

\end{theorem}

The following are useful types of subspaces in symplectic linear algebra.

\begin{definition}\label{annihilator}\index{annihilator}
Let $V$ be a symplectic vector space and $W$ be a linear subspace of $V$. We define its \emph{annihilator}\index{symplectic orthogonal space} as:
\[W^{\circ}:= \{\alpha \in V^{*} \vert\,  \alpha(w)=0, \forall w \in W \}.\]

\end{definition}

\begin{definition}\label{def: symp_orthogonal}
Let $V$ be a symplectic vector space and $W$ be a linear subspace of $V$. We define its \emph{symplectic orthogonal space}\index{symplectic orthogonal space} , by
\[W^{\perp}:= \{v \in V \vert\,  \omega(v,w)=0, \forall w \in W \}.\]

\end{definition}

We can define $\omega^{\# W}\colon V \to W^{*}$ as the restriction of $\omega^{\#}(V)$ to $W$, namely,
\[\omega^{\#}(v)(w):= \omega (v,w), \forall v \in V, w \in W.\]
In this way,
\[W^{\perp}= \ker \omega^{\# W}\] and therefore, we have the induced map $\omega^{\# W}\colon V / W^{\perp} \hookrightarrow W^*$.
In the finite-dimensional setting
\[V / W^{\perp} \simeq W^*\]
and this implies that $\dim W^{\perp}= \dim V -\dim W $.
Denoting by $\omega \vert_ W$ the restriction of $\omega$ to $W$ we get that this induces a bilinear form on $W$ and 
\[(\omega\vert_W)^{\#}= (\omega^{\# W})\vert_W.\] Therefore
\[\ker(\omega \vert_W)^{\#}= W \cap W^{\perp}.\]
Now we define special subspaces of special interest for our purposes.
\begin{definition}\label{def: specialsubsapces}
A subspace $W$ of $V$ is called
\begin{enumerate}
 \item  \textit{Symplectic} \index{symplectic subspace} if $\omega\vert_W$ is a symplectic form or equivalently, 
$W \cap W^{\perp}= \{0\}.$
 \item \textit{Isotropic} \index{isotropic subspace} if $W \subset W^{\perp}$.
\item \textit{Cosiotropic}\index{coisotropic subspace} if $W^{\perp} \subset W$.
\item \textit{Lagrangian}\index{Lagrangiansubspace} if $W^{\perp}=W$. 
\end{enumerate}

\end{definition}

\subsubsection{Symplectic reduction}\label{subsec: symplectic reduction} \index{symplectic reduction}
Let $W$ be a linear subspace of a symplectic vector space $V$. We define
\[ \underline{W}:= W / W \cap W^{\perp},\]
called the \textit{reduction} of $V$. The form $\omega$ induces a symplectic form $\underline{\omega}$ on $\underline{W}$ given by
\[\underline{\omega}([w_1],[w_2]):= \omega(w_1, w_2).\] It can be easily proven that $\underline{\omega}$ is an skew-symmetric, non degenerate form.
The following properties for symplectic reductions and the special subspaces of symplectic spaces hold:
\begin{enumerate}
 \item In the finite dimensional case, if $W$ is a Lagrangian subspace, then $\dim W= \frac 1 2 \dim V$.
 \item If $W$ is coisotropic, 
\[\underline {W}= W / W^{\perp}.\]
\item $W$ is isotropic $\Longleftrightarrow \underline{W}= \{0\}.$ 
\item If $ Z \subset L$, where $L$ is Lagrangian, then $Z$ is isotropic.
\item If $ L \subset Z$, where $L$ is Lagrangian, then $Z$ is coisotropic.
\item If $L \subset Z$, where $L$ is Lagrangian and $Z$ is isotropic, then $Z=L$. That is why a Lagrangian space is also called 
\textit{maximally isotropic}.
\end{enumerate}

\subsubsection{Symplectic manifolds}
Given a smooth manifold $M$, we know the tangent space $T_{p}M$ always has a vector space structure, therefore we can extend the idea of the previous section for a non-linear space (a manifold $M$) defining a symplectic form $\omega_{p}$ on each $T_{p}M$. Let $\omega \in \Omega^{2}(M)$ be a differential 2-form on $M$, then $\omega_{p}\colon T_{p}M \times T_{p}M \mapsto \mathbb{R}$ is bilinear and skew-symmetric on $T_{p}M$ and depending smoothly on $p$. Note that the fact that as $\omega$ is smooth 2-form and non-degenerated, for each function (observable) $f$, there exists a unique vector field $X_f$ so that
$$\iota_{X_f}\omega=df$$

The relation with a Poisson bracket is given by
    \begin{equation*}\label{def: poisson bracket on M}
        \{f,g\} = \omega(X_{f},X_{g})
    \end{equation*}

This definition satisfies all the properties of a Poisson bracket, by definition $\omega$ is skew-symmetric and bilinear, since $X_{f}$ and $X_{g}$ are derivations, the bracket satisfies Leibniz's identity. The geometry of symplectic structures is indeed equivalent to the closedness of the symplectic form. Since $d\omega =0$, then:
\begin{equation*}
    d \omega(X_{f},X_{g},X_{h}) = \{f,\{g,h\}\} + c.p = 0 
\end{equation*}

hence, the Jacobi identity holds. Now, we have the tools to define a Hamiltonian system on any symplectic manifold.

\begin{definition}
    The 2-form $\omega \in \Omega^{2}(M)$ is called $\textbf{symplectic}$ if it is closed ($d\omega = 0$) and $\omega_{p}$ is symplectic for all $p\in T_{p}M$.
    The pair $(M,\omega)$ with $M$ a manifold and $\omega$ a symplectic form is a \textbf{symplectic manifold}. 
\end{definition}

An important fact about symplectic manifolds is that they are always even-dimensional, that is $\text{dim}(M)= \text{dim}(T_{p}M) = 2n$ for $n\in \mathbb{N}$. 

\begin{example}
    If $M$ is any manifold, then the cotangent bundle $T^{*}M$ is endowed with a canonical symplectic form. The symplectic form is given by $\omega=-d\tau_T$ where the 1-form $\tau_T$ in $T^*M$ is known as {\it Liouville form} and is given by  
    $$\tau_T|_{(q,\alpha)}(X)=\alpha_q(d_q\pi X)$$ where $\pi:T^*M\to M$ is the canonical projection of the cotangent fibers.
    This example has a critical application in physics because the configuration space of a physical system is usually defined as a manifold $M$ and its cotangent bundle $T^{*}M$ is known as the phase space.   
$\diamond$\end{example}

\begin{example}\label{ex: canonical form}
    Let's consider $M= \mathbb{R}^{2n}$ with coordinates $(q_{1}, \dots , q_{n}, p_{1},\dots , p_{n})$, then the canonical 2-form:
    \begin{equation}
        \omega_{can} = \sum_{i=1}^{n} dq_{i}\wedge dp_{i} 
    \end{equation}
    is symplectic and canonical isomorphic to $T^*\mathbb{R}^n$.
$\diamond$\end{example}

\begin{definition}
    Let $(M, \omega)$ be a symplectic manifold, and $L \subset M$ a submanifold of $M$, if $T_{x}L$ is a Lagrangian subspace of $T_{x}M$, for any $x\in L$, then $L$ is called a \textbf{Lagrangian submanifold}.
\end{definition}

\begin{example}
    If $S \subset X $ is a submanifold, the "co-normal bundle"
    \begin{equation*}
       N^{*}S = \{ (x, \xi) \vert\ x\in S, \xi \in N_{x}^{*}S\}
    \end{equation*}
    is a Lagrangian submanifold of $T^{*}M$
$\diamond$\end{example}

\begin{example}
   $S^{2}$ is a symplectic manifold and, in fact, the only sphere that admits a symplectic structure. 
$\diamond$\end{example}

\begin{example}
    Given any Lie group $G$, each coadjoint orbit of $G$ on the dual of its Lie algebra is a symplectic manifold. The corresponding symplectic structure  is called the Kirillov-Konstant-Soriau (KKS) symplectic form.
$\diamond$\end{example}
\begin{remark}
    At the beginning of the notes, readers will find a quote from Weinstein about Lagrangian submanifolds. This quote underscores their importance through several applications in geometry and mathematical physics, notably in WKB approximation theory for solving Schr\"odinger's equation (or other linear differential equations). In this context, the phase function is represented as a Lagrangian submanifold in phase space \cite{WBQuanti, maslov1965, hormander1971}. The WKB method and its geometric framework are a crucial bridge between classical and quantum physics.
\end{remark}

Example \ref{ex:standard_symp} and Theorem \ref{thm:Darboux basis} told us that any symplectic vector space could be written in terms of the Darboux basis, this result can be generalized to symplectic manifolds, which means that we can locally see any symplectic form as the form defined in example \ref{ex: canonical form}.

\begin{theorem}[Darboux]\label{thm: Darboux}
    Let $(M, \omega)$ be any $2n-$dimensional symplectic manifold, then for any $q\in M$ there exists a local coordinate chart \[(U, x_{1}, \dots , x_{n}, y_{1}, \dots, y_{n})\] centered at $q$ such that:
    \begin{equation*}
        \omega = \sum_{i=1}^{n}dx_{i}\wedge dy_{i}
    \end{equation*}
\end{theorem}

As a consequence, all symplectic structures are locally identical.

\subsubsection{Hamiltonian systems on symplectic manifolds}
Symplectic geometry provides a natural language to describe classical mechanics, now we can generalize the ideas of section \ref{secham} using symplectic structures.

\begin{definition}[Hamiltonian system]\label{def: ham system}\label{def: ham vec field }
    Let $(M, \omega)$ be a symplectic manifold and let $H\colon M \mapsto \mathbb{R}$ be a smooth function. The unique vector field $X_{H}$ such that $\iota_{X_{H}} \omega = dH$,  is called the \textbf{Hamiltonian vector field} associated to the Hamiltonian function $H$. The triple $(M , \omega, H)$ where $(M,\omega)$ is a symplectic manifold and $H\in C^{\infty}(M)$ is a Hamiltonian function , is called a \textbf{Hamiltonian system}.
\end{definition}
We already gave a general definition for a Hamiltonian vector field in \ref{def: ham poiss}, these two definitions are equivalent since we can define a Poisson bracket for symplectic structures from equation \eqref{def: poisson bracket on M}.

In section \ref{section: integrable systems}, we will study a particular type of Hamiltonian system,which is called integrable, such systems can be solved, provided some conditions.

\subsection{Poisson structures}
Symplectic manifolds are just the non-degenerate case of a Poisson structure, now we will describe the geometry of this structure, the seminal references for these sections are \cite{dufour2005poisson,laurent-gengoux2013poisson,PoissCFM}. Poisson brackets, as defined in section \ref{subsection: poisson bracket}, are bi-derivations, thus they can be represented by {\it bivector fields}. Let $M$ be a smooth manifold, a bivector field on $M$ is simply a section of $\Lambda^2(TM)$. The set of all bivector fields on $M$ will be denoted by $\mathfrak{X}^{2}(M)$. In a local frame, we can write a bivector $\pi \in \mathfrak{X}^2(M)$ as:
\begin{equation}
    \pi(x) = \sum_{i<j}\pi^{ij}(x)\partiald{}{x_{i}}\wedge\partiald{}{x_{j}}
\end{equation}
 We are using the next convention for the wedge product $v \wedge w = v\otimes w -w\otimes v$, which is skew-symmetric. Note that this description implies that the bivector field (pointwise) can be given by a square matrix. The rank of $\pi_x$, the matrix at $x\in M$, is called {\bf the rank} of of $\pi$ at $x$.

The Poisson bracket and the bivector field are related by:
\begin{equation}
    \pi(df,dg) = \{f,g\}.
\label{bivectorPoisson}
\end{equation}

Note that, from skew-symmetry and Leibniz rule, we have a well defined bivector field. However, the Jacobi identity is a missing information in the relation \eqref{bivectorPoisson}. This new data is given in the following theorem:

\begin{theorem}
    The bracket in (3.8) satisfies Jacobi identity if and only if the respective bivector field $\pi$ satisfies 
    \begin{equation*}
        [\pi,\pi] = 0
    \end{equation*}
on the Schouten-Nijenhuis bracket for multi-vector fields.
\end{theorem}

And this result motivate the definition of Poisson manifold:
 \begin{definition}
     The pair $(M,\pi)$, where $M$ is a manifold and $\pi \in \mathfrak{X}^{2}(M)$ satisfying the condition $[\pi,\pi]=0$ is a \textbf{Poisson manifold}.
 \end{definition}

\begin{remark}
    
The Jacobi identity in this context is related to an extension of the Lie bracket for vector fields to the Lie bracket for multivector fields, such an extension is the Schouten-Nijenhuis bracket, and is represented by \cite{moshayedi} :

\begin{multline}
    [X_{1}\wedge \dots \wedge X_{m} , Y_{1}\wedge \dots \wedge Y_{n}] = 
    \sum_{1\leq i \leq m}\sum_{1\leq j \leq n } (-1)^{i+j}[X_{i}, Y_{j}] X_{1}\wedge \dots \wedge X_{i-1}\wedge X_{i+1}\wedge \dots \wedge X_{m} \\
    \wedge Y_{1}\wedge \dots \wedge Y_{j-1}\wedge Y_{j+1}\wedge \dots \wedge Y_{n},
\end{multline}

where $X_{i},Y_{i}$ are vector fields and $[X_{i}, Y_{j}]$ is the usual Lie bracket for vector fields. The Schouten–Nijenhuis bracket satisfies the graded Jacobi identity:
\begin{equation*}
    (-1)^{(|X|-1)(|Z|-1)}[X,[Y,Z]] +  (-1)^{(|Z|-1)(|Y|-1)}[Z,[X,Y]] + (-1)^{(|Y|-1)(|X|-1)} [Y,[Z,X]] = 0
\end{equation*}
In the previous equation, we are using the next notation, $|X|=p$ if $X \in \mathfrak{X}^{p}(M)$

\end{remark}

\subsubsection{Examples of Poisson structures}
In this section, we will show some non-trivial examples of Poisson structures. First, we will construct symplectic foliations for $\mathbb{R}^{3}$, which is a Poisson manifold, and then we will describe the Lie-Poisson structure, from which we can write a Poisson bracket for the dual of a Lie algebra, and with this we will explore some examples. 

\subsection{Lie-Poisson structure}
\label{Lie-Poisson structure}
Another way to construct Poisson structures is through the dual of Lie algebras, some Poisson brackets used in physics are examples of this structure. Let's start with some definitions

\begin{definition}\label{definition: lie algebra}

A {\it Lie algebra} $\mathfrak g$ is a vector space over a field ($\mathbb C$ or $\mathbb R$) provided with a Lie bracket $[\cdot,\cdot]$ that is bilinear, skew-symmetric, and satisfies the Jacobi identity.
$\forall\  x,y,z \in \mathfrak g $ and $\forall\ a,b \in F$:

\end{definition}

\begin{example}
Note that ${\black \mathfrak g }=\mathbb{R}$ is a Lie algebra with $[x,y]=xy - yx = 0$ defined through current multiplication. However, its basis only has one element, $\mathbb{R}={\rm span}(1)$, enough to know that $\mathfrak g$ is abelian: $[1,1] = 0$.
$\diamond$\end{example}

\begin{example}\label{example: so(3)}
    $\mathbb{R}^{3}$ with the vector product is also a Lie algebra which is isomorphic to $\mathfrak{so}(3) = \{A \in GL(3,\mathbb{R}) \colon A = - A^{t}\}$ with the matrix commutator $[A,B]= AB-BA$.
$\diamond$\end{example}

\begin{example}
 $\mathfrak{sl}(n,\mathbb{R}) = \{A \in GL(n,\mathbb{R}) \colon \text{tr}(A)=0\}$ with the matrix commutator. 
$\diamond$\end{example}

Let $\mathfrak{g}$ be any finite-dimensional Lie algebra, then we can give a Poisson structure to the dual of the Lie algebra $\mathfrak{g}^{*}$, this is called the Lie-Poisson bracket and it is given by:
\begin{equation*}
    \{f,g\}(\xi) = \xi([df(\xi),dg(\xi)])
\end{equation*}
Where $f,g \in C^{\infty}(\mathfrak{g}^{*}) $ and $\xi \in \mathfrak{g}^{*}$, the differential $df(\xi) \in T_{\xi}^{*}\mathfrak{g}^{*}$, but since $\mathfrak{g}$ is a finite-dimensional vector space, then $T_{\xi}^{*}\mathfrak{g}^{*} \cong \left( \mathfrak{g}^{*} \right) ^{*} \cong \mathfrak{g}$, so that the bracket $[df(\xi),dg(\xi)]$ is well defined. The term $\xi(\cdot)$ is just the dual pairing between elements of $\mathfrak{g}$ and $\mathfrak{g}^{*}$,that is $\xi(\cdot) = \innerproduct{\xi}{\cdot}$.

\begin{example}
    Let's consider $\mathfrak{g} = \mathfrak{so}(3)$, as we asserted on example \ref{example: so(3)}, $\mathfrak{so}(3)$ can be identified with $\mathbb{R}^{3}$ and the Lie bracket would be the vector product $\times$. The pairing between $\mathbb{R}^{3}$ and $\left(\mathbb{R}^{3}\right)^{*}$ is just the dot product. Then, the Lie-Poisson structure of $\mathfrak{so}(3)^{*}$ is:
    \begin{equation*}
        \{f,g\}(x) = x \cdot \left( \nabla f \times \nabla g \right)
    \end{equation*}
$\diamond$\end{example}
    Note that this structure was used for the Hamiltonian description of the rigid body on example \ref{rigidexample}, which follows naturally since the space of configuration of the rigid body is the rotation group $SO(3)$.
    \begin{example}
        Our construction was meant for finite-dimensional Lie algebras, nevertheless one can find some applications of the Lie Poisson bracket for infinite dimensional spaces. For example, consider an ideal fluid, the configuration space is the group of volume-preserving diffeomorphisms $\mathfrak{D}_{\text{vol}}$, the Lie algebra $\mathfrak{g}$ associated with this group is related to the material representation of the fluid, while $\mathfrak{g}^{*}$ is associated to the spatial description \cite{JMMS}. The dual pairing between is given by integration and the Lie bracket is the Jacobi-Lie bracket. Then the Poisson structure for this system is:
        \begin{equation*}
             \{F,G\}(\textbf{v}) = \int_{D} \textbf{v} \left[\derfunc{F}{\textbf{v}} , \derfunc{G}{\textbf{v}} \right] d^{3}x
        \end{equation*}
        This is the same proposed in example \ref{example: fluid}, but we are explaining its geometric background in this case.
    $\diamond$\end{example}

\subsection{Symplectic leaves}
One important and well-known result of Poisson geometry is that every Poisson structure is locally formed by symplectic manifolds, which we will call symplectic leaves. This statement, referred to as the Weinstein splitting theorem, is a generalization of the Darboux theorem for Poisson manifolds.  
\begin{theorem}[The Weinstein Splitting Theorem]\label{thm: weinstein splitting}
 Around the point $x_{0}$ on a Poisson manifold $(M, \pi)$ whit rank $2k$, there are local coordinates $(q_{1}, \dots , q_{k}, p^{1} , \dots , p^{k} , e_{1} , \dots , e_{l})$, such that:
 \begin{equation*}
     \pi = \sum_{i=1}^{k} \partiald{}{q_{i}}\wedge \partiald{}{p_{i}} + \frac{1}{2} \sum_{1 \leq i < j \leq l} \eta^{ij}(e) \partiald{}{e_{i}} \wedge \partiald{}{e_{j}} 
 \end{equation*}
 given that $ \eta^{ij} (0) = 0 $
\end{theorem}
From this theorem, it is immediate that each leaf of the foliation of a Poisson manifold is symplectic. With this in mind, given a Poisson manifold, we can construct symplectic leaves through its Casimir functions.
\begin{proposition}[Construction of symplectic leaves]\label{proposition: symplectic leaves}
    Let $(M,\pi)$ be a $m$-dimensional Poisson manifold and $p\in M$ regular point (i.e with locally constant rank) with rank $r$. If 
\begin{enumerate}
    \item there are 
$f_1,\dots, f_{m-r}$  Casimir functions 
\item $f_1,\dots, f_{m-r}$ are linearly independent (i.e their differentials, $(df_1)_q , \dots , (df_{n})_{q}$ are linearly independent, for all $q$ in a dense subset of $M$.), 
\end{enumerate}
  then the symplectic leaves are given by connected components of the level sets of $(f_1,\dots, f_{m-r}):M\to \mathbb{R}^{m-r}$. 
\end{proposition}
Using proposition \ref{proposition: symplectic leaves} we will construct some examples of symplectic foliations of $\mathbb{R}^{3}$.
\begin{example}[Foliation by hyperboloids]
Let's consider the next Poisson bivector on $\mathbb{R}^3$:
\begin{equation*}
    \pi = z \partial_{x}\wedge \partial_{y} - 2x \partial_{x}\wedge \partial_{z} + 2y \partial_{y}\wedge \partial_{z} 
\end{equation*}
A Casimir function $f$ must satisfy $X_{f} = \{f , \cdot\} = \pi(df, \cdot) = 0$, that is equivalent to the next system of PDEs:
\begin{eqnarray*}
    z \partiald{f}{x} - 2y \partiald{f}{z}=&0\\
    2x \partiald{f}{z} - z\partiald{f}{y} =& 0 \\
    2y\partiald{f}{y} - 2x \partiald{f}{x} =& 0
\end{eqnarray*}
The solution is the function $f(x,y,z) = 4xy +z^2$. Casimir functions must be constant on the symplectic leaves, so the symplectic leaves are the level sets of $f$, that is:
\begin{eqnarray*}
    4xy +z^2 = c
\end{eqnarray*}
So, we have three options depending on the sign of $c$. Figure \ref{fig:hyperboloids} illustrates some level surfaces of $f$:
\begin{enumerate}
    \item[(i)] $c>0:$ This is a one-sheeted hyperboloid.
    \item[(ii)] $c=0:$ This is a cone. In this case, we have two leaves, the origin and the two connected components of the complement of the origin.
    \item[(iii)] $c<0:$ This is a two-sheeted hyperboloid.
\end{enumerate}
\begin{figure}[htbp]
    \centering
    \begin{subfigure}[b]{0.32\textwidth}
        \centering
        \includegraphics[width=\textwidth]{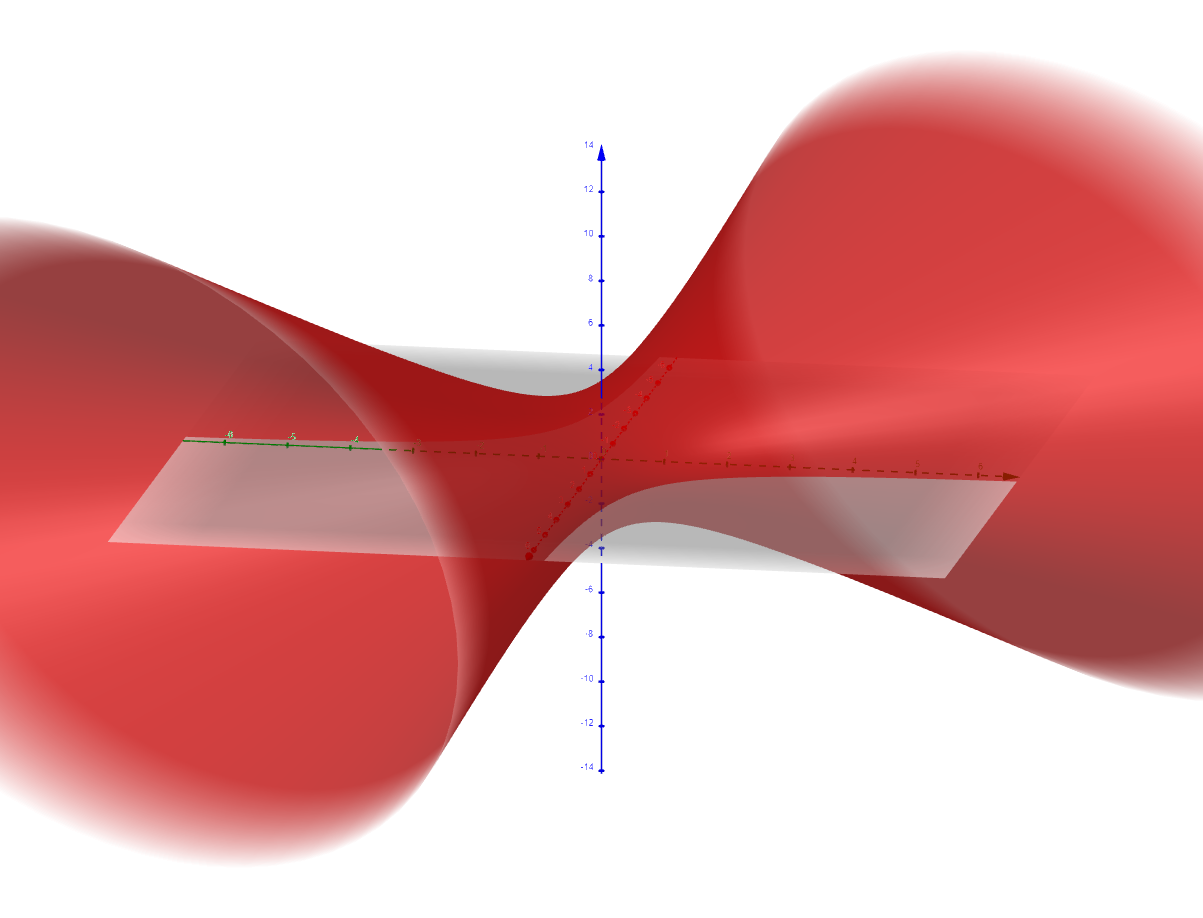}
        \caption{One-sheeted hyperboloid (\( c > 0 \))}
        \label{fig:positive_c}
    \end{subfigure}
    \hfill
    \begin{subfigure}[b]{0.32\textwidth}
        \centering
        \includegraphics[width=\textwidth]{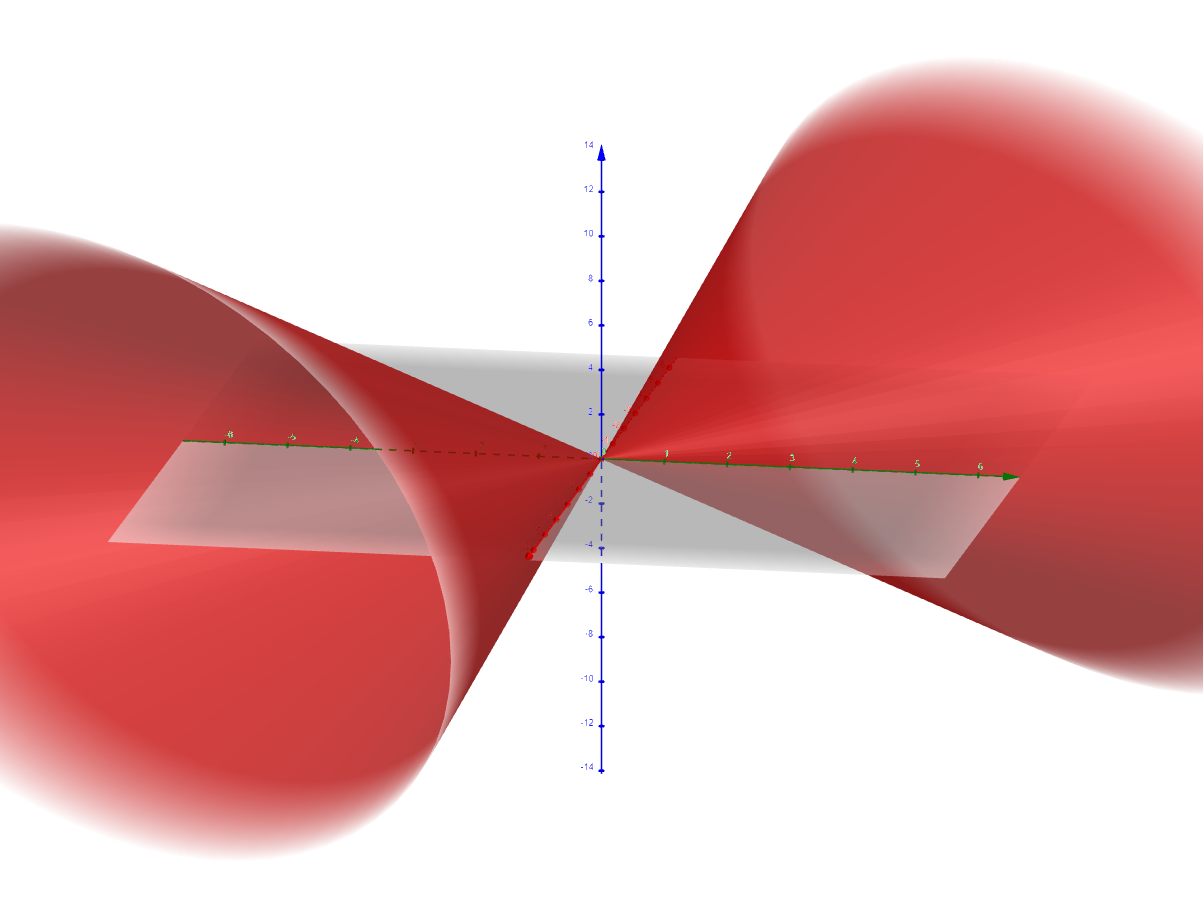}
        \caption{Cone (\( c = 0 \))}
        \label{fig:zero_c}
    \end{subfigure}
    \hfill
    \begin{subfigure}[b]{0.32\textwidth}
        \centering
        \includegraphics[width=\textwidth]{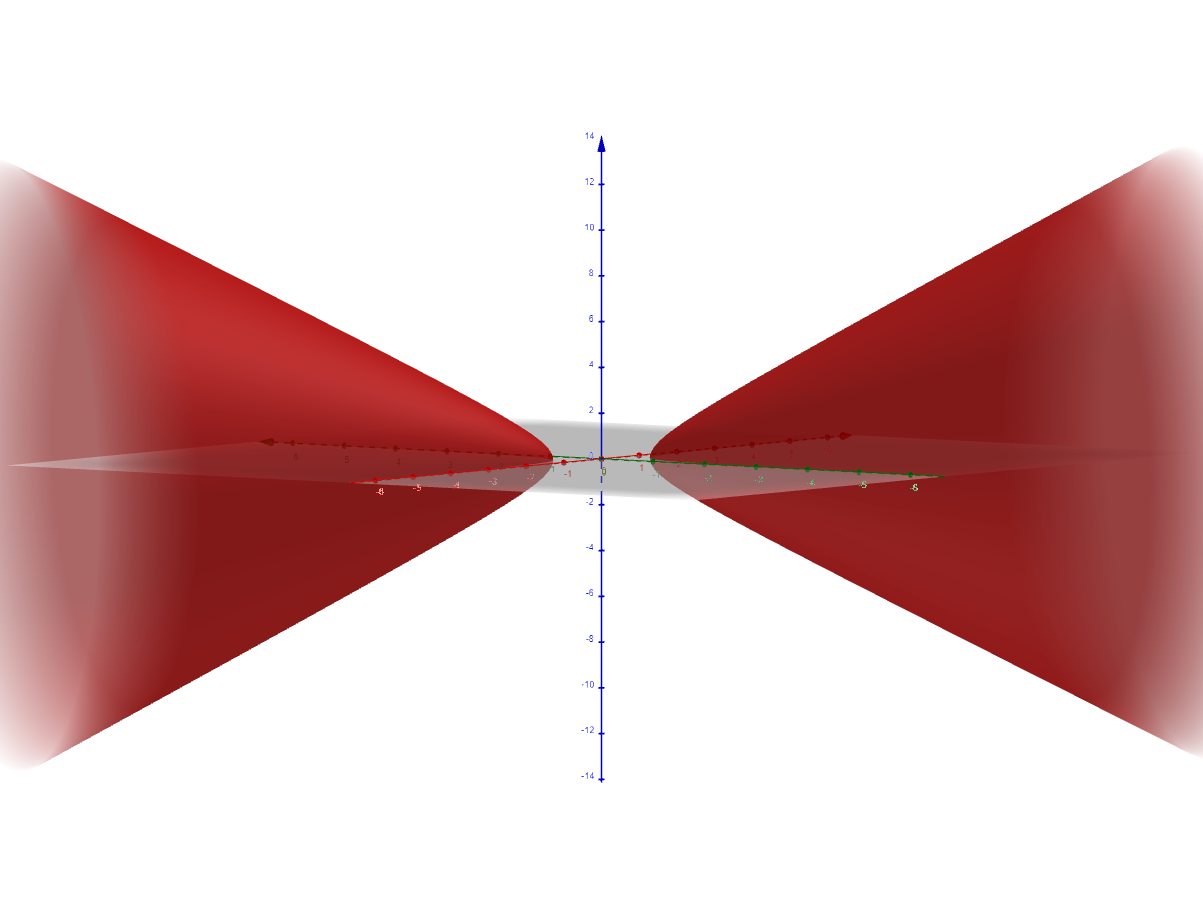}
        \caption{Two-sheeted hyperboloid (\( c < 0 \))}
        \label{fig:negative_c}
    \end{subfigure}
    \caption{Plots of \( c = 4xy + z^2 \) for different values of \( c \).}
    \label{fig:hyperboloids}
\end{figure}

$\diamond$\end{example}
\begin{example}[Foliation by cylinders]
   For this example, let's consider the bivector field on $\mathbb{R}^{3}$, written in cylindrical coordinates, $\pi = r  \partial_{\theta}\wedge \partial_{z}$ and vanishing along $r=0$. Again, any Casimir function $f(r,\theta,z)$ would satisfy the system of PDEs $X_{f}= 0$ , in this case we have:
   \begin{eqnarray*}
       r\partiald{f}{\theta} =& 0 \\
       -r\partiald{f}{z} =& 0
   \end{eqnarray*}
   The solution is a function, which only depends on $r$, namely $f(r,\theta , z) = h(r)$. The level sets of this function are the cylinders $h(r) = c$ , so we have a foliation by cylinders and points along the $z$-axiss. 
$\diamond$\end{example}

\section{And Beyond$\ldots$}

Symplectic and Poisson structures appear in many other contexts apart from those treated in previous sections. In what follows we will briefly explain some of their most relevant applications regarding different areas such as category theory, Lie theory, Dirac structures, hyperbolic geometry, and field theories. We encourage the reader to continue learning and studying the diverse range of applications of symplectic and Poisson geometry.
\\
\subsection{Integrable systems}\label{section: integrable systems}
The solvability of a Hamiltonian system in terms of its integrals of motion was discussed in section \ref{secham}, now we can go deeper on this idea.

\begin{definition}[Integrable system]\label{def: integrable system}
    The Hamiltonian system $(M,\omega, H)$, with $\text{dim}(M)=2n$ is \textbf{integrable} if there exist $n$  independent integrals of motion $f_{1}=H, \dots , f_{n}$, which are pairwise in involution with respect to the Poisson bracket, namely:
    \begin{eqnarray*}
        \{f_{i},f_j\} = \omega(X_{f_{i}},X_{f_{j}}) = 0
    \end{eqnarray*}
    for all $i,j$.
\end{definition}
\begin{theorem}[Arnold \cite{arnold1989mathematical}-Liouville \cite{liouville1855note}]\label{thm: Arnold liouville}
Let $(M,\omega, H)$ be a $2n-$dimensional integrable system with integrals of motion $f_{1}= H, f_{2}, \dots ,f_{n}$, let $c \in \mathbb{R}^{n}$ be a regular value of $f := (H,f_{2}, \dots, f_{n})$ , then $f^{-1}(c)$ is a Lagrangian submanifold, in addition:
\begin{enumerate}
    \item If  the flows of the vector fields $X_{H}, X_{f_{2}}, \dots , X_{f_n}$ starting at a point $q\in f^{-1}(c)$ are complete, then the connected component of $f^{-1}(c)$ containing $q$ is a homogeneous space for $\mathbb{R}^{n}$. With respect to this affine structure, that component has coordinates $\phi_1 , \phi_2, \dots , \phi_n$ known as \textbf{angle coordinates} in which the flows of the vector fields  $X_{H}, \dots , X_{f_n}$ are linear.
    \item There are coordinates $\psi_1, \dots , \psi_n$ known as \textbf{action coordinates}, complementary to the angle coordinates, such that $\psi_{i}$ is a constant of motion for all $i$ and
    \begin{equation*}
        \phi_1 , \phi_2, \dots , \phi_n , \psi_1, \dots , \psi_n
    \end{equation*}
    form a Darboux chart.
\end{enumerate}
\end{theorem}
Thus, any integrable system that satisfies the conditions of the theorem can be solved using action-angle coordinates through quadrature, meaning it can be reduced to a finite number of algebraic operations and integrations.
\begin{example}
    All 2-dimensional phase space, time-independent Hamiltonian systems are integrable . For example, the harmonic oscillator with symplectic form $\omega = dq\wedge dp$ and Hamiltonian function:
    \begin{equation*}
        H(p,q) = \frac{p^2}{2}+\frac{w^2 q^2}{2}
    \end{equation*}
    With constant of motion $f_1 = H$. The Hamiltonian $H(p,q)=E$ level curves are ellipses. This system is integrable, and it can be solved in terms of the action-angle variables:
    \begin{eqnarray*}
        \phi &=& w \tan^{-1}\left( \frac{q}{\sqrt{2\psi-q^{2}}} \right) \\
        \psi &=& \frac{E}{w}
    \end{eqnarray*}
$\diamond$\end{example}

\begin{example}[Spherical pendulum \cite{meinrenken_symplectic}]
Consider a spherical pendulum of length $L=1$ and mass $m=1$. For such a system, the configuration space is $Q = S^{2}$ with coordinates $(\theta,\psi)$, where $\theta \in (0,2\pi)$ and $\psi \in (0,\pi)$. Hence its phase space is $T^{*}S^{2}$, a 4-dimensional symplectic manifold with the symplectic structure $\omega = dp_{\theta}\wedge d\theta + dp_{\psi}\wedge d\psi $. For this system to be integrable, we need at least two constants of motion, one is the Hamiltonian:
\begin{equation*}
    H = \frac{p_{\psi}^{2}}{2}+\frac{p_{\theta}^{2}}{\sin^{2}(\psi)}-g\cos(\psi)
\end{equation*}
Where $p_{\theta} = \sin^{2}(\psi) \dot{\theta}$ and $p_{\psi} = \dot{\psi} $ are the generalized momenta. Observe that $-\partiald{H}{\theta} =\dot{p}_{\theta}=0$, hence $p_{\theta}$ is an integral of motion. In a physical sense,this is related to the rotational symmetry along the rotation axis. Since there are two constants of motion: $H, p_{\theta}$, we conclude that this is an integrable system. 
    
$\diamond$\end{example}

%

The concept of integrable systems is not limited to symplectic manifolds; in fact, it can be extended to integrable systems on Poisson manifolds. Let $(M,\{\cdot , \cdot\})$ be a Poisson manifold, the set of functions $f_{1}, \dots , f_{n}$ is said to be in involution if $\{f_{i},f_{j} \} = 0$ for all $i,j$. Now, we define a Poisson integrable system.

\begin{definition}
    Let $(M,\{ \cdot , \cdot \})$ be a Poisson manifold of rank $2r$ and $F = (f_{1}, f_{2} , \hdots , f_{s} )$, where the functions $f_{1}, \hdots , f_{s}$ is in involution, with $\text{dim}(M)=r+s$. Then $(M,\{ \cdot , \cdot \}, F)$ is an integrable system. The vector fields $X_{f_{i}}$ are called integrable vector fields, and $F$ is known as the momentum map.
\end{definition}

For a comprehensive treatment of this topic, including deeper discussions on integrable systems, and the Poisson version of theorem \ref{thm: Arnold liouville} refer to \cite{Adler_Moerbeke_Vanhaecke, laurentgengoux2008actionanglecoordinatesintegrablesystems}.
\subsection{Intermezzo on extreme-action principle and Lagrangian formulation}
%
In general, the equations of motion of any Hamiltonian system are derivable from an extreme action principle which ensures that the system evolves according to the extremization of an integral of the Hamiltonian function. The action is not only useful to determine the dynamics of a classical system, but also crucial in quantum phenomena theories (to be applied in the calculation of topology invariants), because a particular integral of its exponentiation provides all possible expectation values. Therefore, this interludio may be useful to better understand some of the final topics; see \cite{arnold1989mathematical,JMMS} for much more content on it.
\begin{example}
\label{ex:Actions}
    Considering a one particle system like (\ref{eom}) with time domain in $\mathbb{R}^+$, the {\it action functional} is the integral:
\begin{align}
    S[q,\dot q] = \int_{T \subseteq \mathbb{R^+_0}} L(q,\dot q(p)) 
    ,\quad  \quad  
    L = dt \Bigl[ \dot q(t) p(t) - H(q,p)\Bigr] \in \Omega^1(T)
    \label{Action1dim}
\end{align}
with fixed values $(q,p)|_{\partial T}$ on the boundary of the domain. The integrand is known as the {\it Lagrangian} function, to which we relate the temporal 1-form $L$. When the configuration variables are fields rather than coordinates, like examples \ref{EulerEqn} and \ref{VlasovMaxwell} that have a velocity field $v(x,t):M \subseteq \mathbb{R}^n\times \mathbb{R}^+_0 \to \mathbb{R}^n$, the action is an integral over all the field domain manifold $M$:
\begin{align}
    S[v,\dot v, \partial_x v] = \int_{M} \mathcal 
    {L}(v,\dot v(\eta),\partial_{x}v)  
    ,\quad  \quad  
    \mathcal{L}= dt\wedge dV \left[ \sum_{i=1}^n \dot v_i(x,t) \eta_i(x,t) - \mathcal{H}(v,\eta)\right] \in  \Omega^{n+1}(M)
\label{ActionField}
\end{align}
With $\eta_i = \partial \mathcal{L}/\partial \dot v_i$ the canonical momentum of velocity component $v_i$, both fixed at the boundary $(v,\eta)|_{\partial M}$ and the volume $n$-form $dV = dx_1\wedge\cdots \wedge dx_n$ from Euclidean space $V\subseteq \mathbb{R}^n$. Here $\mathcal{L}$ is the {\it Lagrangian density}, a dim$(M)$-form, and $\mathcal{H}$ is a Hamiltonian density function. 
$\diamond$\end{example}
Note that integrating the Lagrangian density over $V$ provides the Lagrangian 1-form, $\int_V \mathcal{L} = L$; similarly, $\int_V dV \ \mathcal{H} = H$ gives the Hamiltonian function.   
\begin{definition}
    The {\it extreme action principle} states that the actual dynamics of a Hamiltonian system occurs along the trajectories that extremize the action in a variational sense, $\delta S = 0$, and that also satisfy the boundary conditions at $\partial M$. The variational equation gives both the Hamilton equations we have studied and the {\it Euler-Lagrange} equations that constitute an alternative formulation known as {\it Lagrangian}. 
\end{definition}
\begin{example}
    From example \ref{ex:Actions}, since the actions only depend on the variables and their first derivatives, the Euler-Lagrange equations arising from the extremization of (\ref{Action1dim}) and (\ref{ActionField}) are respectively \footnote{Whenever the Lagrangian is a function of higher-order derivatives respect to the coordinates or time, the Euler-Lagrange equations (\ref{EulerLagrangeEqn}) suffer important modifications after the variational extremization.}:
    \begin{equation}
        \frac{d}{dt}\left( \partiald{L}{ \dot q}\right) = \partiald{L}{q}
        \quad {\rm\ and\ }\quad
        \sum_{j,i=1}^n \partial_j\left( \partiald{\mathcal{L}}{(\partial_j v_i)}\right) = \sum_{i=1}^n \partiald{\mathcal{L}}{v_i},
\label{EulerLagrangeEqn}
    \end{equation}
by abusing of the notation and taking $L$
 and $\mathcal{L}$ as integrands without differentials. Therefore, from (\ref{EulerLagrangeEqn}) the Newton equation (\ref{Newton}) follows directly. Furthermore, when writing the Lagrangian in terms of the Hamiltonian by using their Legendre transformation shown in equation (\ref{Action1dim}), and then applying the variational extremization, the Hamilton equations arises. Although being equivalent, the last differential equations are $2n$ of 1st order, in contrast to the equivalent Euler-Lagrange $n$ equations that are of 2nd order.
 $\diamond$\end{example}
%

\subsection{Dirac structures}
As a general structure, that contains symplectic and Poisson structures as particular cases, we found the {\bf Dirac structures}, introduced in \cite{courant1988,courant1990}, and provide a geometric framework that naturally describes Dirac's Constraint Theory \cite{dirac1950}. This makes Dirac structures a great tool to describe the geometry of constrained mechanical systems. 

In brief, a Dirac structure on a vector space \( V \) is defined as a subspace \( D \subset V \oplus V^{*} \) such that a symmetric bilinear form, defined on \( V \oplus V^{*} \), vanishes when restricted to \( D \). This construction extends naturally to a manifold by considering the structure within each tangent space, as follows:

\begin{definition}
    Given a manifold $M$, the {\it generalized tangent bundle}, also referred to as {\it Pontryagin vector bundle}, is the direct sum of its tangent and cotangent bundles, $P = TM \oplus T^*M$. For every $\quad X,Y\in TM,\ \alpha,\beta \in T^*M$, it is endowed with:
\begin{enumerate}
\item A symmetric and bilinear pairing
\begin{equation}
\left< (X,\alpha),(Y,\beta) \right> = \alpha(Y) + \beta(X)
\label{paringTM}    
\end{equation}    
    \item A generalized bracket defined through the Lie bracket and Lie derivatives on $TM$:
    $$ \left[\left[ (X,\alpha),(Y,\beta) \right]\right] = (\ [X,Y]\ ,\ \mathcal{L}_X \beta - \mathcal{L}_Y \alpha 
 + \frac{1}{2}( d\alpha(Y) + d\beta(X))\ ).$$
 \end{enumerate}
\end{definition}
\begin{definition}
    A {\it Dirac structure} on $M$ is a Lagrangian submanifold  $D \subset P$, this is $D =D^\perp$ respect to the pairing (\ref{paringTM}), such that the generalized bracket between elements of smooth differential sections $\Gamma(D)$, lies on them: $[\![ \Gamma(D),\Gamma(D) ]\!] \subset \Gamma(D)$.
\end{definition}
\begin{example}
    Let $M$ be a manifold, $\omega$ be a 2-form on $M$, and $\Delta_{M}$ a distribution on $M$, that is $\Delta_{M}$ is a vector subbundle of $TM$. Now, the next subbundle $D_{\Delta_{M}} \subset TM\oplus T^{*}M$:
    \begin{multline*}
        D_{\Delta_{M}}(m) := \{(v_{m},\alpha_{m}) \in T_{m}M\times T_{m}^{*}M : v_{m} \in \Delta_{M}(m) \text{ and } \\ \innerproduct{\alpha_{m}}{y_{m}} = \omega(m)(v_{m},y_{m})
        \text{ for all } y_{m} \in \Delta_{M}(m) \}
    \end{multline*}
    is a Dirac structure on $M$.
    \label{example: dirac1}
$\diamond$\end{example}
Moreover, we can describe dynamical systems on Dirac structures, for example, if $D$ is a Dirac Structure on a manifold $M$, and $f : M \mapsto \mathbb{R}$ a smooth function, the dynamical system for a curve $\sigma(t) \in M$ is given by the equations $(\dot{\sigma}(t), df(\sigma(t))) \in D(\sigma(t))$. Constructions as the one on example \ref{example: dirac1} are very relevant to describing non-holonomic \footnote{Non-holonomic constraints are those that can't be expressed in terms of the positions and the time, they depend on quantities such as the velocity, for this reason, they are not integrable. A non-holonomic constraint can be written as $f(q,\dot{q},t) = 0$. For example, a ball of radius $R$ rolling on a plane that is rotating with a varying angular velocity along the vertical axis, is a non-holonomic system.} and implicit Lagrangian systems \cite{yoshimura2006}, or even for extending the classical Lagrange-d'Alembert variational principle to more complex variational principles as the Lagrange-d'Alembert-Pontryagin or the Hamilton-Pontryagin principle as in \cite{MR2265469}.    

To demonstrate the power of Dirac structures in modeling constrained systems, we will now focus on recent advancements in the geometry of non-equilibrium thermodynamics, as discussed in the works of Gay-Balmaz and Yoshimura \cite{GayYoshimura2018, GayYoshimura2020}. We will highlight some key findings from these papers, but a more comprehensive discussion and additional constructions can be found in the original references.

First, we limit the description to the case of \textit{simple thermodynamical systems}, which are macroscopic systems that can be completely described in terms of a finite set of mechanical variables $(x_{1}, \hdots , x_{n})$ and one thermal variable, that can always be chosen as the entropy. We start with the thermodynamic configuration space $\mathrm{Q} = Q \times \mathbb{R}$, where $Q$ is the mechanical configuration space with coordinates $q$, and the entropy $S \in \mathbb{R}$. The Lagrangian $L(q,\dot{q},S) : TQ\times \mathbb{R} \mapsto \mathbb{R}$, and let's assume the presence of a friction force and an external force $F^{\text{fr}}, F^{\text{ext}} : TQ\times \mathbb{R} \mapsto T^{*}Q$. An extension of the classical variational principle to a non-holonomic variational principle for non-equilibrium thermodynamics was proposed in \cite{GayYoshimura2017}. A curve $(q(t), S(t)) \in \mathrm{Q}$, is a solution for such variational formulation if it satisfies:
\begin{equation}
\delta\int_{t_1}^{t_2}L(q,\dot{q},S)dt + \int_{t_1}^{t_2}\innerproduct{F^{\text{ext}}(q,\dot{q},S)}{\delta q}dt = 0
\end{equation}

\begin{equation}
    \partiald{L}{S}(q,\dot{q},S) \delta S = \innerproduct{F^{\text{fr}}(q,\dot{q},S)}{\delta q}
    \label{eqn: variational constraint}
\end{equation}

\begin{equation}
    \partiald{L}{S}(q,\dot{q},S) \dot{S}= \innerproduct{F^{\text{fr}}(q,\dot{q},S)}{\dot{q}}
    \label{eqn: phenomenological constraint}
\end{equation}
We know from statistical physics that $ \partiald{L}{S}(q,\dot{q},S)= -T$, with $T$ the temperature of the system, which implies that $ \partiald{L}{S}(q,\dot{q},S)<0$. Every variation $\delta q(t)$ and $\delta S(t)$ are restricted to the variational constraint \ref{eqn: variational constraint}. The non-linear non-holonomic constraint \ref{eqn: phenomenological constraint} is known as a phenomenological constraint, because it is determined by experimental aspects of the friction force. The variational constraint \ref{eqn: variational constraint}, gives rise to the next submanifold of $T \mathrm{Q}\times_{\mathrm{Q}} T \mathrm{Q} $:
\begin{equation}
    C_{V} = \left\{(q,S,v,W,\delta q, \delta S) \in T \mathrm{Q}\times_{\mathrm{Q}} T \mathrm{Q} \colon \partiald{L}{S}(q,v,S)\delta S = \innerproduct{F^{\text{fr}}(q,v,S)}{\delta q} \right\}
    \label{eqn: var cons}
\end{equation}
On the other hand, from the phenomenological constraint \ref{eqn: phenomenological constraint} we get the following constraint:
\begin{equation}
     C_{K} = \left\{(q,S,v,W) \in T \mathrm{Q} \colon \partiald{L}{S}(q,v,S)W = \innerproduct{F^{\text{fr}}(q,v,S)}{v} \right\}
\end{equation}

We know from example \ref{example: dirac1} that we can construct a Dirac structure on the Pontryagin bundle $\mathrm{P}$ through a distribution. Furthermore, in a general case, we can obtain that distribution using the variational constraint $C_{V}$ as follows. Let $(q,p,v) \in \mathrm{P}$, then we can induce a distribution, which is locally defined as:
\begin{equation}
    \Delta_{\mathrm{P}}(q,p,v) = \left \{(q,v,p,\delta q, \delta v, \delta p) \in T_{(q,v,p) }\mathrm{P}: (q, \delta q) \in C_{V}(q,v) \right \}
    \label{eqn: induced distribution}
\end{equation}
Following \ref{eqn: induced distribution}, let $x = (q,S,v,W,p,\Gamma) \in \mathrm{P} = T \mathrm{Q}\oplus T^{*} \mathrm{Q} $ then, for the variational constraint $\ref{eqn: var cons}$, we get the next distribution
\begin{equation*}
    \Delta_{\mathrm{P}}(x) = \left\{(x,\delta x) \in T\mathrm{P
    } : \partiald{L}{S}(q,v,S)\delta S = \innerproduct{F^{\text{fr}}(q,v,S)}{\delta q}   \right\}
\end{equation*}
As in example \ref{example: dirac1}, the Dirac structure on $\mathrm{P}$ induced from $\Delta_{P}(x)$ is given by:
\begin{equation*}
   D_{\Delta_{\mathrm{P}}}(x) = \{(v_{x},\alpha_{x}) \in T_{x}\mathrm{P}\times T_{x}^{*}\mathrm{P} : v_{x} \in \Delta_{\mathrm{P}}(x) \text{ and } \innerproduct{\alpha_{x}}{y_{x}} = \omega_{\mathrm{P}}(x)(v_{x},y_{x}) \text{ for all } y_{m} \in \Delta_{\mathrm{P}}(x) \}
\end{equation*}
In this case, $\omega_{P}$ is the presymplectic form \footnote{A presymplectic form $\alpha$ is a closed differential 2-form of. If $\alpha$ is also non-degenerate, it defines a symplectic structure.} induced from the canonical symplectic form $\Omega_{T^{*} \mathrm{Q}}$ on $T^{*} \mathrm{Q}$, that is $\omega_{\mathrm{P}} = \pi^{*}_{(\mathrm{P},T^{*}\mathrm{Q})}\Omega_{T^{*} \mathrm{Q}}$, where $\pi_{(\mathrm{P},T^{*}\mathrm{Q})}\colon \mathrm{P} \mapsto T^{*}\mathrm{Q} $. Locally, $\omega_{\mathrm{P}}$ can be computed as:
\begin{equation}
    \omega_{\mathrm{P}}(q,S,v,W,p,\Gamma) = dq_{i}\wedge dp_{i} + dS \wedge d\Gamma
\end{equation}
Finally, a curve $x(t) = (q(t),S(t),v(t),W(t),p(t),\Gamma(t)) \in \mathrm{P}$ satisfies the Dirac dynamical system:
\begin{equation}
    ((x,\dot{x}), d\varepsilon(x)) \in D_{\Delta_{\mathrm{P}}}(x)
    \label{eqn: dirac dynamical}
\end{equation}
Where $\varepsilon: \mathrm{P} \mapsto \mathbb{R}$ is known as the \textit{generalized energy} on the Pontryagin bundle, and for a Lagrangian $L(q,v)$ is provided by:
\begin{equation*}
    \varepsilon(q,v,p) := \innerproduct{p}{v}-L(q,v)
\end{equation*}
The condition \ref{eqn: dirac dynamical} is equivalent to the evolution equations for the thermodynamics of non-equilibrium simple systems, namely:
\begin{eqnarray}
    \frac{d}{dt}\partiald{L}{\dot{q}} - \partiald{L}{q} &=& F^{\text{ext}}(q,\dot{q},S) + F^{\text{fr}}(q,\dot{q},S) \\
    \partiald{L}{S}\dot{S} &=& \innerproduct{F^{\text{fr}}(q,\dot{q},S)}{\dot{q}}
\end{eqnarray}

Though the previous examples just consider constraints systems, this is not the only application of such structures. We also can find them on generalized geometry \cite{gualtieri2011generalized} or Manin pairs \cite{bursztyn2008courant}, among others.

\subsection{The symplectic category}
We proceed to define and discuss some version of the category of symplectic manifolds. In particular, we show schematically how to construct a category of symplectic manifolds such that the composition of its morphisms is always well defined. We refer the reader to \cite{SympWeinstein} for a concise discussion and references. 
\begin{definition}
The category $ {\mathbf L_{ \mathbf{ Symp}} }$ of symplectic vector spaces is defined to have ${\rm Obj}({\mathbf L_{\mathbf{Symp} }}) \ni (V,\omega)$ whose morphisms are {\it linear symplectomorphisms} ${\rm Mor}({\mathbf L_{\mathbf{Symp} }}) \ni f : (V,\omega_V) \to (W,\omega_W)$ from definition \ref{symplectomorphism}, all maps for which the symplectic structure is preserved.
\end{definition}
\begin{definition}
    The category of symplectic manifolds $\mathbf{Symp }$ is composed by manifolds with a non-degenerated and closed symplectic 2-form, ${\rm Obj}( \mathbf{ Symp}) \ni (M,\omega)$, and the morphisms between them ${\rm Mor}( \mathbf{ Symp}) : (M,\omega_V) \to (N,\omega_N)$ are symplectomorphisms.
\end{definition}
In particular, note that if $f:(V,\omega)\to(V,\omega)$ is a linear symplectomorphism its graph is Lagrangian ${\rm gr}(f) = {\rm gr}(f)^\perp \subseteq (V,\omega) \times (V,-\omega)$. This suggests the possibility of extending ${\rm Mor}({\mathbf L_{\mathbf{Symp} }})$ to include all Lagrangian subspaces.
%
\\~\\
Equivalently, we may extend $\mathbf{ Symp}$ to a new category $\mathbf{ Symp}^\mathbf{ext}$ whose objects continue to be symplectic manifolds, but whose morphisms are Lagrangian submanifolds. Initially, $\mathbf{Symp}^\mathbf{ext}$ was not a category because the compositions among their morphisms was not completely closed. One way to solve this issue is via the Wehrheim-Woodward construction (WW)  introduced in \cite{WW} and summarized in what follows.
\begin{definition}
    The category $\mathbf{Span}$, has all sets as objects and its morphisms are triples of sets related by functions. Be $X,Y,Z\in {\rm Obj}(\mathbf{Span})$ three sets, an associated morphism $F\in {\rm Mor}(\mathbf{Span})Y$ is $X \leftarrow Z \rightarrow Y$.
\end{definition}
\begin{definition}
    The {\it category of relations} between sets, denoted $\mathbf{Rel}$, has all sets as objects, while its morphisms are all the subsets of the cartesian product of the two objects. For $X,Y\in {\rm Obj}(\mathbf{Rel})$ two sets, a morphism $Y\xrightarrow{F} X$ is $F\subseteq X\times Y$ and it is called a {\it relation}. 
\end{definition}
\begin{definition}
\label{def: RelComposition}
    Given two relations $f,g\in {\rm Mor}(\mathbf{Rel})$ such that $f \subseteq X \times Y$ and $g \subseteq Y \times Z$, their {\it composition} is the relation:
    \begin{equation*}
        X\times Z \supseteq f \circ g = \{ (x,z)\ |\ (x,y)\in f,\ (y,z)\in g 
     {\rm\ with\ a\ common\ } y\in Y \},
    \end{equation*}
    While their {\it fiber product} is:
    \begin{equation*}
    f \times_{Y} g 
    = f \times g \ \cap \  X\times \Delta_{Y} \times Z,
    \end{equation*}
    with the diagonal set $ Y\times Y :=  \Delta_{Y} = \{ (y,y) | y\in Y \} $ that plays the role of the identity morphism $Y\xrightarrow{\Delta_{Y} } Y$.
\end{definition}
\begin{definition}
    A {\it canonical relation} is a Lagrangian and smooth submanifold $L\subseteq M\times N^-$, such that $N,M\in {\rm Obj}(\mathbf{ Sym})$ and we denote $N^- = (N,-\omega_N)$.
\end{definition}
Here smoothness refers to closeness, then a canonical relation must be a closed submanifold.
Note that symplectic manifolds and canonical relations are also objects and relations from $\mathbf{ Rel}$, respectively. 
\begin{definition}
    Two submanifolds $F,G\subset M$ are {\it transversal} if the direct sum of their tangent spaces at every point is $ T_p F + T_p G = T_p M$, and $\forall p\in M$.
\end{definition}
\begin{definition}
    The pair of canonical relations $(f,g)$ such that $f \subseteq X \times Y^-$ and $g \subseteq Y^- \times Z$ are said to be:
    \begin{enumerate}
        \item {\it Transversal} if $f \times g$ is transversal to the product $X\times \Delta_{Y^-}\times Z$.
        \item {\it Strongly transversal}, $f \pitchfork g$, if additionally, the projection of the fibre product $f \times_{(Y^-)} g$ into $X \times Z$ is an embedding onto a closed manifold, which is $f\circ g$ with the composition definition (\ref{def: RelComposition}).
\end{enumerate}
\end{definition}
Transversality is the condition we require to ensure the suitable composition between canonical relations. In particular, if the pair $(f,g)$ is transversal, then $f \times_{(Y^-)} g$ is a manifold; but if $f \pitchfork g$, then $f\circ g$ is a closed manifold and also a canonical relation. To determine the transversality of canonical relations we apply the next criterion.
\begin{theorem}
    The pair of canonical relations $f \subseteq X \times Y^-$ and $g \subseteq Y^- \times Z$ is transversal if and only if their tangent bundles $Tf \subseteq TX \times TY^-$ and $Tg \subseteq TY^- \times TZ$ are such that the domain of $Tf$ is transverse to the range of $Tg$.
\end{theorem}
The WW construction is now schematically explained introducing a couple of new concepts.
\begin{definition}
    A finite {\it canonical path} $f=(f_1,\ldots,f_i)$ is a sequence of composable canonical relations, in the left direction, whose source and target are those from $f_1$ and $f_i$ respectively. Their composition is defined to be concatenation whenever it is admissible: $f \circ g = (f_1,\ldots,f_i,g_1,\ldots,g_j)$ only if the target of $g_1$ is the source of $f_i$. And the identity path related to every $X\in \mathbf{ Symp}$ is a sequence of identity relations $(\Delta_X,\ldots,\Delta_X)$.\end{definition}
\begin{definition}
\label{WWcanonicalRelation}
    Two canonical paths $f$ and $g$ are related, $f\sim g$, if one is obtained from the other by replacing some of the successive entries $(f_i,f_{i+1})$ with the entry $f_{i}\circ f_{i+1}$, only if $f_{i} \pitchfork f_{i+1}$.
\end{definition}
Therefore, each pair of canonical relations $(f)$ and $(k)$ whose composition $h\circ k$ is a canonical relation can be collapsed into the same entry $(h\circ k)$. Also, the equivalence class of every canonical path is represented by the canonical relation obtained by collapsing all the composable entries.
The quotient of the set of canonical paths by the previous equivalence relation lets to the following category.
\begin{definition}
    The morphisms of the {\it Wehrheim-Woodward category} ${\rm WW}(\mathbf{Symp^\mathbf{ext} })$ are all the equivalence classes $[(f_1,\ldots ,f_r)]$ of canonical paths under $\sim$ defined in \ref{WWcanonicalRelation}, and their composition is the equivalence class of the concatenated paths $[(f_1,\ldots ,f_r,g_1,\ldots,g_r)]$. Their objects are the same than ${\rm Obj(\mathbf{Symp})}$.  
\end{definition}
\begin{example}
    Every canonical relation $f\subseteq X \times Y^-$ is the canonical path of one entry $(f)$. Representatives of its equivalence class are formed by all paths with tails of identity relations, $[(f)]\ni(\Delta_X,\ldots,\Delta_X,f,\Delta_{Y^-},\ldots,\Delta_{Y^-})$ , considering that they compose suitably with $f$. Because $\Delta_X \pitchfork f$ and $f\pitchfork\Delta_{Y^-}$, we are sure that $\Delta_X \circ f$ 
 and $f\circ 
 \Delta_{Y^-}$ are closed Lagrangian submanifolds.
$\diamond$\end{example}
\begin{definition}
    The $composition functor$ between categories $c: {\rm WW}(\mathbf{Symp^\mathbf{ext} }))\to \mathbf{ Rel}$ acts as identity on symplectic manifolds, the objects, and maps each equivalence class of canonical paths, the morphisms, to the canonical relation obtained by composition of its entries:
    $$
    {\rm Mor}({\rm WW}(\mathbf{ Symp}))\ni [(f_1,\ldots ,f_r)] \mapsto f_1\circ f_2 \circ \cdots \circ f_r \in {\rm Mor(\mathbf{Rel})}
    $$
\end{definition}
The composition is now closed between any pair of canonical paths, and its analog in  $\mathbf{Symp}$ or $\mathbf{Rel}$ is defined by using the composition functor. However, a relevant question still to be answered explicitly is how to obtain invariants of symplectic manifolds via ${\rm WW}(\mathbf{ Symp})$. Recent results on the versions of the symplectic category, the WW construction and its connections with topological field theory can be found in \cite{hcc, ckm, CMS}
\subsection{Mapping class group and hyperbolic geometry}
As we will see, this group is important in the study of hyperbolic geometry and 3-manifolds, also as an invariant of Chern-Simons theory covered below. To get more insight we refer the reader to the lecture notes \cite{Minsky2013}, or \cite{FarbMargMCG} for a more detailed discussion.
\begin{definition}
    The {\it mapping class group} of a manifold $M$ is defined as the quotient:
    $${\rm MCG}(M) = {\rm Aut}(M) / {\rm Aut}^0(M) $$
    where ${\rm Aut}(M)$ are the automorphisms on $M$, a group under the operation of composition, and ${\rm Aut}^0(M)$ are all of them that belong to the connected component where the identity automorphism lies. 
\end{definition}
In the context of surfaces $S$, ${\rm Aut}(S)$ may be replaced by the group of homeomorphisms ${\rm Hom}(S)$. Then, the identity's connected component is equivalent to the identity isotopy class: automorphisms for which exist a homotopy, being also a homeomorphism, to the identity.
\\~\\
For every surface $S_g$ of genus $g\geq 1$ its MCG$(S_g)$ has a symplectic representation on linear integer symplectomorphisms, Sp$(2g,\mathbb{Z})$, whose matrix representation has integer entries; each isotopy class of homomorphisms is then represented by one symplectic matrix. This is because the surface first homology group, $H_1(S_g,\mathbb{Z})$, can be extended to a symplectic vector space whose symplectic form is constructed as follows.
\\
\begin{definition}
    The {\it algebraic intersection number} of two transverse, oriented, simple closed curves $\gamma$ and $\eta$ in an oriented surface $S$, denoted as $\hat i (\gamma,\eta)$, is the sum of their intersections indices. Where an intersection has index $+1$ if it has the same orientation as $S$, and $-1$ in other cases.
\end{definition}
The intersection number can be reduced to that between the homology classes $[\gamma],[\eta] \in H_1(S_g,\mathbb{Z})$, any pair of closed curves from each class will have the same intersection number. For any $\alpha \in [\gamma]$ and $\beta\in [\eta]$ then $\hat i (\alpha,\beta)=\hat i ([\gamma],[\eta])$. Considering that
$$
    \hat i : H_1(S_g,\mathbb{Z}) \times H_1(S_g,\mathbb{Z}) \to \mathbb{Z}
$$
is bilinear, skew-symmetric, and nondegenerate, when extending it to $\mathbb{R}$, the first homology group is naturally endowed with a symplectic structure, $(H_1(S_g,\mathbb{R}) , \hat i)$. Moreover, its automorphisms are symplectomorphisms respecting $\hat i$, and they have a representation in symplectic matrices Aut$( H_1(S_g,\mathbb{Z})) \to {\rm Sp}(2g,\mathbb{Z})$. In consequence, the natural representation of the mapping class group on the first homology group MCG$(S_g)\to  {\rm Aut}(H_1(S_g,\mathbb{Z}))$ is then traduced to a representation on Sp$(2g,\mathbb{Z})$.
\\
\\
The Darboux basis of the extension $(H_1(S_g,\mathbb{R}) , \hat i)$ is the geometrical symplectic basis $\{[\gamma_1],..,[\gamma_g],[\eta_1],..,[\eta_g]\}$ \footnote{For this homology classes basis the algebraic and geometric intersection numbers are the same. The last one is simply the minimum possible number of intersections between representatives of each class.}, where $\gamma_n$ and $\eta_n$ are curves surrounding the  $n$-th hole and $1\leq n \leq g$. From the related dual basis, the 2-form related to the algebraic intersection number is written:
$$
\omega_{can} = \hat i = \sum_{n=1}^g  [\gamma_n]^*\wedge [\eta_n]^*
$$
%
%
\subsection{Lie Theory}

This theory studies Lie groups, $G$, generated by Lie algebras, {\black $\mathfrak g$} already introduced in section \ref{Lie-Poisson structure} to construct Poisson structures. We suggest, for example, \cite{Kirillov2008,Warner1983,weyl1946classical} as a first introductory reading and deeper ones respectively. The generation of the Lie group occurs in such a way that $G$ plays the role of the global object, {\black a manifold that respects the structure of a group}, while {\black $\mathfrak g$} is its local version, {\black the tangent space of that manifold}. For this reason, the passing from Lie groups to Lie algebras is called {\it differentiation}, whereas the inverse process is an {\it integration}, both being functors between the respective categories. 
\\
A Lie algebra may be thought of as simpler than a Poisson algebra since it also has a bilinear operation without a Leibnitz rule.
Because it is a vector space $\mathfrak g$ with dimension $n$, it has a basis, $\{x_1,..,x_n\}$ composed of {\it generators}, we only must know the action of $[\cdot,\cdot]$ on it.
Whenever the Lie bracket is zero $\forall x,y \in \mathfrak g$, which means that the algebra is commutative, we call it {\it abelian}.
\begin{definition}
A {\it Lie group} $G$ is both a manifold and a group whose multiplication $\mu : G \times G \to G,\ (g,h)\mapsto g*h$, and inverse maps $\iota : G \to G, \ g\mapsto g^{-1}$, are smooth: $\mu, \iota \in C^\infty(G)$.
\end{definition}
\begin{example}
    The circle $G = S^1 \subset \mathbb{C}$ is a Lie group with rotation as the operation. Because $S^1 = \{z\in \mathbb C\ ;\ |z|=1 \}$, each element is $z = e^{i\phi}$ with $\phi\in [0,2\pi)$, their product $z_1\cdot z_2 = e^{i(\phi_1+\phi_2)}$ is a rotation from one on the other, and the identity is $1 = e^{i0}$. {\black Justly, the Lie algebra of $S^1$ is ${\mathfrak s}^1=\mathbb{R}$, that is, it is generated infinitesimally by the the complex exponentiation of the the algebra: $z = e^{i\phi}$.} 
$\diamond$\end{example}
\begin{example}
Considering its applications to the knots invariants and particle physics symmetries, we discuss the Lie Group of special unitary matrices, first of degree 2:
\begin{align}
    SU(2) &= \{ U \in M_{2\times 2}(\mathbb{C});\ U^{-1} = U^\dagger,\ \det(U)=1  \}.
\end{align}
With matrix multiplication as the group operation. $SU(2)$ has a dimension of $3$ given $\dim(SU(N)) = N^2-1$, and we may reproduce it by using the fundamental representation of matrices $2\times2$ or representations on higher order matrices. The associated Lie algebra  ${\mathfrak{su}(2)}$ must be a 3-dimensional vector space of $2\times 2$ hermitian and traceless matrices with the Lie bracket built from matrix multiplication \footnote{Another convention is taking skew-hermitian generators, $i\sigma/2$ instead $i\sigma$, in which case the Lie bracket is $[x_i, x_j] =  \sum_{k=1}^3\epsilon_{ijk} x_k$}:
\begin{align}
    {\black {\mathfrak{su}(2)} } &= \{ X \in M_{2\times 2}(\mathbb{C});\ X^\dagger = X,\ {\rm Tr}(X)=0  \} = {\rm span}(\{x_1, x_2, x_3\})
    \label{su(2)algebra}
\end{align}
\begin{align}
    [x_i, x_j] &= x_ix_j - x_j x_i =  i \sum_{k=1}^3\epsilon_{ijk} x_k 
    \label{su(2)bracket}
\end{align}
Where $\epsilon_{ijk}$, the totally-skew-symmetric tensor, are the structure constants of $SU(2)$. When working in the fundamental representation the Lie algebra is spanned out by the Pauli matrices,
\begin{align}
    \mathfrak{su}(2) = {\rm span}
    \Bigl( \Bigl\{ 
    \sigma_1=
    \left(    \begin{array}{cc}
        0 & 1 \\
        1 & 0 \\
    \end{array}\right)
    ,
    \sigma_2 =
    \left(   \begin{array}{cc}
        0 & -i \\
        i & 0 \\
    \end{array}\right)
    ,
    \sigma_3=
   \left(     \begin{array}{cc}
        1 & 0 \\
        0 & -1 \\
    \end{array}\right)
    \Bigr\} \Bigr),
\label{su(2)span}
\end{align}
and by convention it is chosen as generators $x_i = \sigma_i/2$. Then, each vector in $\mathfrak{su}(2)$ is $\sum_{k=1}^3 a_k \sigma_k /2$ with $a_k \in \mathbb{C}$
,  and every element of $SU(2)$ is the imaginary exponentiation of an element in $\mathfrak{su}(2)$:
\begin{align}
    SU(2)\ni U = e^{i \sum_{k=1}^3 u_k \sigma_k/2} = \sum_{p=0}^\infty \frac{(i)^p}{p!} \left( \sum_{k=1}^3 u_k \sigma_k/2 \right)^p, \ u_k \in \mathbb{C}.
\end{align}
$\diamond$\end{example}
\begin{example}
In general, the $SU(N)$ Lie group of unitary matrices with determinant one have $N^2-1$ generators that span its Lie algebra ${\mathfrak su}(N)$, composed of hermitian and traceless matrices. The fundamental representation is naturally on matrices $N\times N$\footnote{In the standard model of particles, the electromagnetic and weak interactions are invariant under transformations from $U(1) \otimes SU(2)$, while the weak and strong interactions governing the nuclear dynamics are invariant respect to $SU(3)$. The complete symmetry is $U(1) \otimes SU(2) \otimes SU(3)$ to be related to the famous N\"oether's theorem that states that any symmetry of a Hamiltonian system is equivalent to a conservation law. In this case the conserved quantities are weak hypercharge and weak isospin, that become electric charge conservation after a process called symmetry breaking. Lastly, the color charge from quarks that build up hadrons like protons and neutrons, but that is not physical since has infinite versions because is not gauge invariant \cite{Schwartz2013}.}
$\diamond$\end{example}
Notice that that the series expansion of the exponential until certain finite order shows the infinitesimal character of the Lie algebra respect to the full exponential, infinite series, that represents a Lie Group element. Furthermore, the exponentiation of the Lie algebra implies a continuous generation of all group elements around the Lie Group identity. The next theorem proceeds naturally.
\begin{theorem}[Lie's 3rd theorem]
For every Lie algebra {\black $\mathfrak g$} there is a Lie Group $G$ such that the Lie algebra is the tangent space of that Lie group {\black at the group identity $e$}: $\mathfrak g\ {\black \cong T_eG}$. 
\end{theorem}
However, the {\black {\it integration} } of $\mathfrak g$ has not a unique $G$, same as $\mathbb{R}^n$ is the tangent plane of many $n$-dimensional manifolds.
 Employing category theory we say that between the category of Lie algebras $\rm{\bm Alg}_{\rm Lie}$ (vector spaces with a Lie bracket) and the category of Lie groups $\rm{\bm{ Grp}}_{\rm Lie}$ there exists exist two functors, an {\it integration}$: \mathfrak{g}  \mapsto G$ and a {\it differentiation}$:G \mapsto \mathfrak g$ called the Lie functor, one inverse of the other.
%
\subsubsection{Algebroids and groupoids}
Natural parametrized versions of Lie algebras and Lie groups are Lie algebroids and Lie groupoids respectively, then they serve as generalization to them. We would like to refer the reader to the brief lecture notes \cite{Bursztyn2023,Meinrenken2017} or well to the book \cite{Mackenzie2005a}.
\begin{definition}
A {\it Lie algebroid} is a vector bundle $A$ on a manifold $M$, $A \to M$, and an anchor map $\rho: A \to TM$ between vector bundles, such that the space of sections $\Gamma(A) = \{ \sigma: M \to A \}$, is a Lie algebra whose bracket satisfies a Leibniz rule:
$$ [x, fy] = (\rho(x) f) y + f [x,y],$$
for every pair $x,y\in \Gamma(A)$ and a smooth function $f \in C^\infty(M)$.
\end{definition}
\begin{example}
All Lie algebras forced to be tangent spaces over a point $\mathfrak g \to \{ {\rm pt} \}$ are Lie algebroids. The anchor map sends all elements to the point $\rho:\mathfrak g \to \{ {\rm pt} \}$, the only function is the identity on $\{pt\}$, and the bracket becomes that one from $\mathfrak g$.
$\diamond$\end{example}
\begin{example}
The tangent bundle of any manifold {\black $TM\to M$} is a Lie algebroid whose anchor is the identity $\rho = I: TM \to TM$, and the Lie bracket is the subtraction of the commutation between vector fields on $TM$.
\label{LieAlgebroidTMtoTM}
$\diamond$\end{example}
Other relevant examples are Lie algebra bundles $\mathfrak{g} \times M \to M$, and cotangent bundles {\black $T^*M\to M$ of a Poisson manifold $(M,\pi)$ with $\rho = \pi^\sharp$}.
%
\begin{definition}
A {\it Lie groupoid} is a groupoid whose elements belong to a certain manifold $M$ while its arrows/morphisms are in another manifold $G$. That is a category where all arrows from $G$ have an inverse and act on $M$. Moreover, the following data regarding the operation on $G$.
\begin{enumerate}
 \item Two structure maps, source and target, $s,t:G \rightrightarrows M$ from the manifold of morphisms $G$ to the  manifold of objects $M$, such that each arrow $g\in G$ relates its source with its target $t(g) \xleftarrow{g} s(g)$.
 \item An associative operation $\mu$ on $ G\times G\supseteq G_2 = \{(g_1 , g_2) \in  G\times G \ |\ s(g_1) = t(g_2)\}$ implying that $\mu: G_2 \to  G, \ (g_1,g_2) \mapsto g_1*g_2
$ is restricted to composable arrows.
\item An inverse map $\iota: G \to G,\ g \mapsto g^{-1}$, providing the inverse arrow respect to $*$, $\mu(g,\iota(g)) = g*g^{-1}= 1_{s(g)}$ and $ \mu(\iota(g) , g) = g^{-1}*g = 1_{t(g)} $. 
\item An unit map $\varepsilon:M \to G$, providing the identity arrow for every $x\in M$.  It is unital in the sense that given $\varepsilon(s(g)) = 1_{s(g)}$ and $\varepsilon(t(g)) = 1_{t(g)}$, then $\mu(g , \varepsilon(s(g)) )= g = \mu( \varepsilon(t(g)) , g ) $. Also rewritten as $g*1_{s(g)} = 1_{t(g)}*g = g$.
\end{enumerate}
\end{definition}
\begin{example}
By definition, a groupoid may have infinite identity elements respect to $\mu$, $\{1_{x}\}_{x\in M}$. However, in the simplest example, any Lie group with trivial structure maps onto a point, $G \rightrightarrows \{{\rm pt} \}$, all its elements $g \in G$ have the same source and target, $\mu$ is the group multiplication, and there is just one identity, $1_{\rm pt}$. 
$\diamond$\end{example}
\begin{example}
The fundamental groupoid $\pi_1(M) \rightrightarrows M$, where $\pi_1(M) = \bigl\{[\gamma:[0,1] \to M]\bigr\}$ are all free homotopy-path equivalence classes on $M$ \footnote{$[\gamma]\in \pi_1(M)$ is the equivalence class of paths on $M$ for which exist a homotopy to path $\gamma$ and whose initial and final points are $\gamma(0)$ and $\gamma(1)$; that is, all continuously deformable paths between two ''fixed'' points}, is justly the {\it integration} of the Lie algebroid $TM \to M$ of example \ref{LieAlgebroidTMtoTM}. Here $s([\gamma]) = \gamma(0)$ and $t([\gamma]) = \gamma(1)$, the operation between classes is the concatenation of paths, and the inverse $[\gamma]^{-1}$ is the class of contrary orientated paths.
\label{fundamentalGroupoid}
$\diamond$\end{example}
Other examples of Lie groupoid are the following: 
manifolds over themselves $M \rightrightarrows M$ where all arrows and structure maps are identities relating each point with itself. Pair groupoids $M \times M \rightrightarrows M$ that coincides with the fundamental groupoid of $M$ (example \ref{fundamentalGroupoid}) if $M$ simply connected, because every pair of points are related without referring to a homotopy-path class defined with them. Lastly, Lie group bundles $G\times M \rightrightarrows M$ are natural examples of Lie groupoids.
%
\\~\\
In general one cannot ensure the existence of an {\black {\it integration}} from any Lie algebroid $A$ to a Lie groupoid $G$; thus, there is not an analog of Lie's 3rd theorem that guarantees the existence of a Lie groupoid $G$ such that $A \cong TG|_{\varepsilon(M)}$. The topological obstructions concerning this nonexistence were first presented in \cite{INTLIE}. However, the existence and generalizations of Lie's 3rd theorem  remains to be an open area of research, including infinite dimensional algebroids and {\black groupoids}, higher category theory, for instance 2-{\black groupoids}, and the local and global symmetries of field theories \footnote{A global symmetry is an invariance of the theory under a global transformation, which does not depend on the coordinates in the domain. A local symmetry is regarding a local transformation, those that indeed depend on the point in the domain.}. 

A first connection of this theory with Poisson geometry is, as we state above, $(T^*M,\pi^*,[\cdot,\cdot]_{SN})$ is a Lie algebroid and (when it is integrable) it is integrable to symplectic groupoids, which are Lie groupoids $G$ with a symplectic form that is compatible with the groupoid structure maps. Such integration result also works for other geometries. Interesting connections with more theory of symplectic and Poisson geometry can be considered, for example geometric quantization, reduction of structure, symmetries and dynamics of mechanical systems, among others.
\subsection{TQFT and knot invariants}
A knot $K$ is an embedding of $S^1$ into a three-manifold $(M,g)$. We restrict ourselves to knots in the 3-sphere, $K\subset S^3$, on which a Chern-Simons allows us to recover some knot invariants like Jones and Kauffman polynomials. This is done using the path integral formalism from Quantum Field Theory (QFT) that relies on a special weighted average over all possible configurations of the field, which is a connection on $S^3$. Find references about QFT as it is introduced below, while for TQFT we recommend the lecture notes \cite{labastida1997,freed1992,TQFTSchwarz}.   
\begin{definition}
    A {\it Topological Quantum Field Theory} (TQFT) of Schwarz-type is a QFT whose action, functional of the fields on a manifold $(M,g)$, and all observables expectation values, to define later, do not depend on the metric $g$.
\end{definition}
\begin{example}
    A Chern-Simons theory (CST) on a 3-dimensional manifold $M$ is a Schwarz TQFT whose field is a connection 1-form, $A\in \Omega^1(M)$, defined on a Lie group, $G$, bundle and taking values on its non-abelian Lie algebra $\mathfrak{g}$. Its action is:
\begin{align}
    S_{CH}[A] = \frac{k}{2\pi} \int_{M}  {\rm Tr } \left( \frac{1}{2} A \wedge dA + \frac{1}{3} A\wedge A\wedge A \right),
\label{CSaction}
\end{align}
$\diamond$\end{example}
a bulk integral of a 3-form, with $k$ an integer parameter of the theory, and where the trace is taken on the matrix representation of $\mathfrak g$.
\\
Writing the connection on the cotangent space leads to $A = \sum_{k=1}^3 A_k dx^k$ with every $A_k \in \mathfrak g$ and ${\rm span}(\{dx^i\}_{i=1}^3)=T_p^*M $ for every $p\in M$. Moreover, if the algebra is $\mathfrak{su}(2)$, then $A_k = \sum^3_{b=1}  A_k^b \sigma_b/2$ in the fundamental representation (\ref{su(2)span})\footnote{The Lagrangian form, integrand in (\ref{CSaction}), may be written as: \\ $
\mathcal{L}_{CS} =\frac{1}{2}\sum_{i<j<k}^3 \left( \frac{1}{2}\sum_{b=1}^3A_i^b(\partial_j A_k^b - \partial_k A_j^b ) +  \frac{i}{3}  \sum_{a,b,c=1}^3 A_i^a A_j^b A_k^c \epsilon_{bca} \right) dx^i\wedge dx^j\wedge dx^k,$ \\ by using (\ref{su(2)bracket}) and ${\rm Tr}(\sigma_a\sigma_b)/4 = \delta_{ab}I/2$, satisfied by all fundamental generators of $\mathfrak{su}(N)$.}. Notice that we have never used the metric inside the SC action, its SC form is completely metric-free, and then it encodes topological information of the manifold $M$.
\begin{definition}
    A {\it gauge transformation} is a mapping of the connection fields through elements of its algebra: $A \mapsto A + d\lambda + [\lambda , A]$, where $\lambda \in \mathfrak{g}$ is a 0-form, a function, and $[\cdot,\cdot]$ is the Lie bracket. This is nothing but the infinitesimal version of a general transformation based on the Lie group elements $g = \exp( \lambda) \in G$, with which $A \mapsto g^{-1}Ag + g^{-1} dg$\footnote{Written in the Lie algebra $\mathfrak{su}(2)$ components, it maps the connection to: \\ $\frac{1}{2}\sum_{j,b}\bigl[ (A_j^b + \partial_j\lambda^b)\sigma_b + i \sum_{a,c}\lambda^a A_j^b \epsilon_{abc}\sigma_c\bigr]dx^j$.}.
\end{definition}
The action (\ref{CSaction}) is gauge invariant regarding transformations on the connected component of the identity gauge transformation, while regarding others, it is invariant modulo a multiple of $2\pi$ thanks to $k\in\mathbb{Z}$.
\begin{definition}
    We call -{\it observable}- any quantity that is invariant under gauge transformations\footnote{The word {\it observable} comes from the fact that measurable quantities of nature are gauge invariant in physical theories.}, and denote $\mathcal{A}\subseteq \Omega^1(M)$ the space of 1-forms, $A \in \mathfrak{g}$, that admit a gauge transformation.
\end{definition}
In the classic domain, when extremizing the action (\ref{CSaction}) respect to $A$ using (\ref{EulerLagrangeEqn}), the resultant equation of motion is that one from flats connections, those with zero curvature: $0=\delta_A S_{CS} = dA + A\wedge A := F$. With this usual curvature definition from a connection, $F$ is a 2-form belonging to $\mathfrak{g}$; in particular $F =  \sum_{i,j,b} F^b_{ij}\ \sigma_b/2\  dx^i\wedge dx^j$ with components $F^b_{ij} = \partial_i A_j^b - \partial_j A_i^b + [A_i,A_j]$.
\\~\\
Some knot invariants arise from calculating the expectation value of its Wilson loop, which is the trace of the holonomy of the connection $A$ along the knot. Such interpretation was shown to work in \cite{Witt} employing a 3D CST with $M=S^3$ and algebras $\mathfrak{su}(2)$ and $ \mathfrak{so}(2)$ respectively for Jones and Kauffman polynomials. 
The holonomy is the solution to the parallel transport differential equation, solved by an {\black infinite-iterative} procedure resulting in an exponential\footnote{The parallel transport differential equation respect to a path $\gamma(s):[0,1]\to M$ becomes| $U_\gamma^\prime(s) = i \gamma^\prime(s)\cdot A\ U_\gamma(s)$ with $x\cdot A = \sum_i^3 x_iA_i$.}. 
\begin{definition}
    The parallel transport operator along a knot $K$ is its {\it holonomy}, a path-ordered and imaginary exponentiation of the line integral of the connection $A$ along $K$:
    \begin{align}
    U_K[A] = \mathcal{P}
    \exp 
    \left(i \oint_K A
    \right)
    =
    \sum_{n=0}^\infty \frac{i^n}{n!} \left(
    \oint_{K} \sum_{k_n=1}^3 A_{k_n} dx^{k_n} \cdots \oint_{K} \sum_{{k_1}=1}^3 A_{k_1} dx^{k_1}
    \right)
\end{align}
Where $\mathcal P$ is the path-ordered product operator that orders from left to right according to the partition $K = \cup_{i=1}^n \gamma_i$, from the last path $k_n$ to the first one $k_1$, and is indispensable since $\mathfrak g$ is non-abelian.
\end{definition}
\begin{example}
    If the connection is $A \in\mathfrak{su}(N)$, for every knot the holonomy is $U_K[A]\in SU(N)$ because it is the exponentiation of an element in the Lie algebra. 
$\diamond$\end{example}
\begin{definition}
    The {\it Wilson loop} of a knot $K$ is the trace of its holonomy, $W_K[A] = {\rm Tr}_\rho(U_K[A])$, taken on a matrix representation of the Lie group, $G \to \rho$. Because we restrict ourselves to the fundamental representation, we will not use the subscript $\rho$.
\end{definition} 
Now we introduce the path integral formalism. It relies on a partition or generating function analog to the probability theory one, but in an infinite dimensional case. It integrates over all possible variable values (connection fields) a weight functional somewhat similar to the probability density (exponential of the action (\ref{CSaction})). To see a formal approach and principal references see {\cite{PathIntegral}}, while to study the approach used in particle physics see {\cite{Srednicki2007}} for advanced and self-contained chapters and {\cite{Schwartz2013}} for a detailed explanation from basic topics. 
\begin{definition}
Given the action functional $S_{CH}[A]$ from which the classical equations of motion of field $A$ are derived by extremization, the {\it partition function} is the path integral:
\begin{align}
    Z = \int_{\mathcal{A}\subseteq \Omega^1(M)} \mathcal D[A]\ e^{i S_{CS}[A]}, 
    \quad \quad
    \mathcal D[A] \approx \prod_{k=1}^3 \prod_{b=1}^3 \prod_{x\in M} dA^{b}_k(x) 
\label{pathintegral}
\end{align}
Where the integration occurs on all admissible configurations of the field. From a heuristic construction, the path integral measure is the limit of the product of all possible connection measures on $M$, along the three directions, and on the Lie algebra.  
\end{definition}
In general, the solution of a path integral comes in analogy to the Gaussian integrals of coupled variables in the finite-dimensional case: a function of the determinant of the matrix of couplings, that in the infinite-dimensional case becomes a functional determinant of the differential operator inside the action.
\begin{definition}
    The {\it expectation value} from any observable quantity, that is gauge invariant, is the path integral calculated by plugging it as part of the integrand in (\ref{pathintegral}), and normalizing the result with the partition function.
\end{definition}
We recall the similitude with expectation values in finite-dimension probability theory.
\begin{example}
    The natural quantities being topological and gauge invariant are the Wilson loops. The expectation value of a Wilson loop of a knot $K$ is:
\begin{align}
    \bigl< W_K \bigr >  = \frac{1}{Z} \int \mathcal D [A]\ {\rm Tr}\left(U_K[A]\right)\ e^{i S_{CS}[A]}.
\label{WilsonLoop}
\end{align}
$\diamond$\end{example}
\begin{example}
In the general case, given a link formed out by non-intersecting knots, $L = \bigcup_{i=1}^n K_i$, the Wilson link is the product of the Wilson loops from all its knots. Its expectation value is $\bigl< W_L \bigr >  = \bigl< W_{K_1} \cdots W_{K_n} \bigr > $.
$\diamond$\end{example}
Finally, when setting $M=S^3$, $A\in \mathfrak{su}(2)$ and then $U_K[A]\in SU(2)$, the expectation values (\ref{WilsonLoop}) and $\bigl< W_L \bigr >$ give respectively the Jones polynomial of the knot $K$ and link $L$:
\begin{equation}
    V_L(q^{1/2}(k,2)) = \left< W_L\right>(k) \quad {\rm with}\quad 
    q(k,N) := \exp\left({\frac{2\pi i}{N+k} }\right).
\end{equation}
Where $k$ is the integer parameter of the action (\ref{CSaction}) and $N$, set to $2$, is the dimension of the fundamental representation of the Lie group $SU(2)$. Apart from the links invariants we may recover by means of expectation values, by itself the partition function $Z$ defined in \ref{pathintegral} is an invariant of the complete manifold $M$.

\subsection{The Poisson Sigma Model}
To conclude these notes, we offer a quick overview of another example of how symplectic and Poisson geometry appear naturally in physics.

The Poisson sigma model (PSM) is a topological field theory that corresponds to the following data:
\begin{enumerate}
\item A compact surface $\Sigma$, possibly with boundary, called the \textit{source}. \index{source}
\item A finite dimensional Poisson manifold $(M,\Pi)$, called the \textit{target}. \index{target}
\end{enumerate}



The space of fields for this theory is denoted with $\Phi$ and corresponds to the space of vector bundle morphisms of class $\mathcal C^{k+1}$ between $T\Sigma$ and $T^*M$.
This space can be parametrized by a pair $(X, \eta)$, where \[X\in \mathcal C^{k+1}(\Sigma, M)\]  and \[\eta \in 
\Gamma^k(\Sigma, T^*\Sigma \otimes X^*T^*M),\] and $k \in \{0,\,1,\, \cdots \, , \infty\}$ denotes the regularity type of the map that we choose to work with.\\
On $\Phi$, the following first order action  is defined:
\[S(X,\eta):= \int_{\Sigma} \langle \eta,\, dX\rangle+  \frac 1 2 \langle \eta, \, (\Pi^{\#}\circ X) \eta \rangle,\]
where,

 \begin{eqnarray}
 \Pi^{\#}&\colon& T^*M\to TM\\
  & &\psi \mapsto \Pi(\psi, \cdot).
 \end{eqnarray}
 

In the action functional, $dX$ and $\eta$ are regarded as elements in $\Omega^1(\Sigma, X^*(TM))$, $\Omega^1(\Sigma, X^*(T^*M))$, respectively and  $\langle \,,\, \rangle $ is the pairing between $\Omega^1(\Sigma,  X^*(TM))$ and $\Omega^1(\Sigma,  X^*(T^*M))$ induced by the 
natural pairing between $T_xM$ and $T_x^*M$, for all $x \in M$.

Cattaneo and Felder \cite{cattaneo-felder:poissonsigma} proved that the reduced phase space of the Poisson Sigma Model, if smooth, when the surface is a disk, has the structure of a symplectic groupoid. For more recent work on the connections bewteen PSM and symplectic groupoids, see \cite{Cattaneo2014RelationalSG, Contreras2015}.

This model also has significant importance for deformation quantization \cite{kontsevich2003deformation}. Namely, the perturbative expansion through Feynman path integral 
for PSM, in the case that $\Sigma$ is a disk, gives rise to Kontsevich's star product formula, i.e. the semiclassical expansion of the path integral
\begin{equation}
\int_{X(r)=x}f(X(p)) g (X(q))\exp{(\frac{i}{\hbar} S(X,\eta)) dX d\eta}
\end{equation}
around the critical point $X(u)\equiv x, \eta \equiv 0$, where $p,q$ and $r$ 
are three distinct points of $\partial \Sigma$, corresponds to the star product $f\star g (x)$.

\bibliographystyle{abbrv}
\bibliography{BMBibTeX.bib}

\end{document}